\documentclass[11pt]{amsart}

\usepackage[T1]{fontenc}
\usepackage[utf8]{inputenc}
\usepackage{lmodern}
\usepackage{microtype}
\usepackage{amsmath,amssymb,amsthm,mathtools,mathrsfs}
\usepackage{enumitem}
\usepackage{aliascnt}
\usepackage{needspace}
\usepackage[colorlinks=true,urlcolor=blue,linkcolor=black,citecolor=blue]{hyperref}
\usepackage[nameinlink,noabbrev]{cleveref}
\allowdisplaybreaks
\numberwithin{equation}{section}

\newtheorem{theorem}{Theorem}[section]
\newaliascnt{proposition}{theorem}
\newtheorem{proposition}[proposition]{Proposition}
\aliascntresetthe{proposition}
\newaliascnt{lemma}{theorem}
\newtheorem{lemma}[lemma]{Lemma}
\aliascntresetthe{lemma}
\newaliascnt{corollary}{theorem}
\newtheorem{corollary}[corollary]{Corollary}
\aliascntresetthe{corollary}
\theoremstyle{definition}
\newaliascnt{definition}{theorem}
\newtheorem{definition}[definition]{Definition}
\aliascntresetthe{definition}
\theoremstyle{remark}
\newaliascnt{remark}{theorem}
\newtheorem{remark}[remark]{Remark}
\aliascntresetthe{remark}
\crefname{theorem}{theorem}{theorems}
\Crefname{theorem}{Theorem}{Theorems}
\crefname{proposition}{proposition}{propositions}
\Crefname{proposition}{Proposition}{Propositions}
\crefname{lemma}{lemma}{lemmas}
\Crefname{lemma}{Lemma}{Lemmas}
\crefname{corollary}{corollary}{corollaries}
\Crefname{corollary}{Corollary}{Corollaries}
\crefname{definition}{definition}{definitions}
\Crefname{definition}{Definition}{Definitions}
\crefname{remark}{remark}{remarks}
\Crefname{remark}{Remark}{Remarks}

\newcommand{\R}{\mathbb R}
\newcommand{\Q}{\mathbb Q}
\newcommand{\Z}{\mathbb Z}
\newcommand{\N}{\mathbb N}
\newcommand{\C}{\mathbb C}
\newcommand{\Lio}{\mathscr L}
\newcommand{\Lfield}{\mathbb L}
\newcommand{\Eset}{\mathcal E}
\newcommand{\Ggap}{\mathcal G}
\newcommand{\abs}[1]{\left|#1\right|}
\newcommand{\norm}[1]{\left\lVert#1\right\rVert}

\newcommand{\den}{\operatorname{den}}
\newcommand{\height}{\mathrm h}
\newcommand{\cF}{c_{\mathrm F}}
\newcommand{\Clow}{C_{\mathrm{low}}}

\title[Mahler's Problem II. Rational Maps]
{Mahler's Problem on Liouville Numbers
II. Rational Maps}
\author{Diego Marques}
\thanks{Supported by CNPq, Grant No.~304467/2023-5.}
\address{Departamento de Matem\'atica, Universidade de Bras\'ilia,
Bras\'ilia, 70910-900, Brazil}
\email{diego@mat.unb.br}
\date{}
\hypersetup{
 pdftitle={Mahler's Problem on Liouville Numbers II. Rational Maps},
 pdfauthor={Diego Marques},
 pdfsubject={An arithmetic characterization of real rational maps preserving Liouville numbers},
 pdfkeywords={Liouville numbers, rational functions, arithmetic rigidity, coefficient field}}

\begin{document}

\begin{abstract}
In 1984, Mahler asked which analytic functions map the set $\Lio$
of Liouville numbers into itself. In the first paper of this
series, we proved that every real-analytic solution on an open
interval is the restriction of a real rational function. We now
identify a canonical subfield $\Lfield$ of $\R$ such that, for
every nonconstant $R\in\R(X)$ and every nonempty open interval
$I$ in its domain, $R(I\cap\Lio)\subseteq\Lio$ if and only if
$R$ is defined over $\Lfield$ and
$R^{-1}(\Q)\cap I\cap\Lio=\varnothing$. The weaker inclusion $R(I\cap\Lio)\subseteq\Q\cup\Lio$
holds if and only if $R$ is defined over $\Lfield$.
If $R$ is not defined over $\Lfield$, every such interval
contains a Cantor set of Liouville numbers whose images under
$R$ are irrational and have irrationality exponents bounded
by a constant depending only on $R$. We also construct continuum
many positive strong Liouville numbers $\sigma$ for which
$\sigma X^m$ preserves $\Lio$ for every $m\geq2$, with all
nonzero rational fibers consisting of $U_m$-numbers.
\end{abstract}

\maketitle

\setcounter{tocdepth}{1}
\tableofcontents

\section{Introduction}
\label{sec:introduction}

A real number $\xi$ is a \emph{Liouville number} if, for every integer
$N\geq1$, there exist infinitely many rational numbers $p/q$, with
$q\geq2$, such that
\begin{equation}\label{eq:liouville}
 0<\abs{\xi-\frac pq}<q^{-N}.
\end{equation}
We denote the set of these numbers by $\Lio$. A classical theorem of
Maillet \cite{Maillet1906} states that every nonconstant rational
function with rational coefficients maps $\Lio$ into itself.
Motivated by this theorem, Mahler posed the following question
in 1984 \cite[\S~1, p.~101]{Mahler1984}.
\begin{quote}
\emph{Which analytic functions $f(z)$ have the property that if $\xi$
is any Liouville number, then so is $f(\xi)$? In particular, are there
entire transcendental functions with this property?}
\end{quote}
We say that a function has \emph{Maillet's property} if it preserves
every Liouville number in its domain. Throughout this paper, $D_R$ denotes the real domain of a rational
function $R$, obtained by excluding its real poles. Preservation
is understood only at points of $D_R$.

Mahler's question was subsequently recorded by Waldschmidt
\cite[Question~3.27]{Waldschmidt2004} and Bugeaud
\cite[Problem~44]{Bugeaud2004}.
Earlier work constructed transcendental analytic or entire
functions mapping distinguished uncountable subsets of
$\Lio$ into $\Lio$
\cite{LelisMarquesRamirez2017,MarquesMoreira2015,MarquesSchleischitz2016}; the behavior of Liouville numbers
under rational powers was also studied in
\cite{MarquesSchleischitz2016}.
More recently, the finite-smoothness and $C^\infty$
settings have been shown to admit considerably greater
flexibility \cite{LelisMoreiraSilva2026,MarquesSmooth2026};
see the introduction of Part~I \cite{MarquesI}
for further background.

Mahler's question has both an existence aspect and a classification
aspect. The local rigidity theorem of Part~I
\cite[Corollary~1.2]{MarquesI} shows that every
real-analytic function on an interval with Maillet's property is the
restriction of a real rational function. In particular, an entire
function with this property must be a polynomial, giving a negative
answer to the question about transcendental entire functions.
This determines the analytic form of a solution, but leaves its
arithmetic structure to be understood:
\emph{which rational functions with real coefficients have Maillet's
property?}

The answer cannot be obtained simply by restricting the coefficients
to $\Q$. In degree one, a real M\"obius transformation has
Maillet's property on a nonempty open interval precisely when
it belongs to $\Q(X)$, as we prove in \Cref{cor:mobius}.
In every higher degree, however, there are preserving monomials
with transcendental coefficients. Indeed, \Cref{thm:strong-nonlinear}
produces continuum many positive strong Liouville numbers $\sigma$
for which $\sigma\xi^m\in\Lio$ for every $\xi\in\Lio$ and
$m\geq2$. Their nonzero rational fibers have a precise description
in terms of Mahler's classes $U_m$. A classification must explain
both the degree-one rigidity and this higher-degree flexibility.

We give such a classification in terms of a coefficient field and
the arithmetic of rational fibers. The decisive step is to separate
the two ways in which preservation can fail: a Liouville number may have an irrational image of finite irrationality exponent, or it may
have a rational image. Put
\[
 \Eset:=\Q\cup\Lio.
\]
We prove that, for a nonconstant real rational function $R$, the
relaxed inclusion
\begin{equation}\label{eq:weak-intro}
 R(\Lio\cap D_R)\subseteq\Eset
\end{equation}
is equivalent to $R$ being defined over a specific subfield
$\Lfield\subset\R$. Maillet's property is then characterized by this
coefficient condition together with the exclusion of Liouville
points from the rational fibers. Combined with the analytic rigidity
theorem of \cite[Corollary~1.2]{MarquesI}, this yields a necessary and sufficient
coefficient-and-fiber criterion for the nonconstant real-analytic
solutions to Mahler's problem.

\subsection{The coefficient field}

The relevant arithmetic property concerns not only the accuracy of
rational approximations, but also the scales at which they are
available. Ordinary Liouville approximation allows exceptionally good
approximations to occur at arbitrarily sparse denominator scales.
For the coefficient field, we require approximations of each fixed
order to be available throughout all sufficiently large scales,
within a polynomial denominator bound.

\begin{definition}\label{def:field}
Let $\Lfield$ be the set of real numbers $\alpha$ with the following
property: for every integer $N\geq1$, there exist integers $D\geq1$
and $H_0\geq2$ such that, for every integer $H\geq H_0$, there are
$p\in\Z$ and $q\in\Z_{>0}$ satisfying
\begin{equation}\label{eq:scale-property}
 q\leq H^D,
 \qquad
 \abs{\alpha-\frac pq}\leq(Hq)^{-N}.
\end{equation}
\end{definition}

The exponent $D$ may depend on the prescribed approximation order
$N$. This flexibility allows approximations of different coefficients
to be synchronized with a sufficiently accurate rational approximation to a given $\xi\in\Lio$. Exact rational approximations are permitted, so
$\Q\subseteq\Lfield$.

The definition is purely Diophantine and makes no reference to the
action of a function on $\Lio$. We prove that $\Lfield$ is a field and
that
\begin{equation}\label{eq:field-inclusions-intro}
 \Q\subsetneq\Lfield\subsetneq\Eset.
\end{equation}
Its role as a coefficient field is established by the following
classification.

\subsection{Classification and local escape}

For a nonconstant $R\in\R(X)$, write
\begin{equation}\label{eq:normalized-pair}
 R=P/Q,\qquad P,Q\in\R[X],\qquad \gcd(P,Q)=1,
 \qquad Q\text{ monic}.
\end{equation}
This reduced, normalized pair is unique, and
$D_R=\{x\in\R:Q(x)\ne0\}$. With this normalization, $R\in\Lfield(X)$ is equivalent
to $P,Q\in\Lfield[X]$.

\Needspace{10\baselineskip}
\begin{theorem}[Rational classification]\label{thm:classification}
Let $R\in\R(X)$ be nonconstant. The following conditions are equivalent.
\begin{enumerate}[label=\textup{(\roman*)},leftmargin=2.7em]
\item $R\in\Lfield(X)$.
\item $R(\Lio\cap D_R)\subseteq\Eset$.
\item There is a nonempty open interval $V\subseteq D_R$ such that
      $R(\Lio\cap V)\subseteq\Eset$.
\end{enumerate}
\end{theorem}

Thus the inclusion $R(\Lio\cap V)\subseteq\Eset$ on a single
nonempty open interval $V\subseteq D_R$ forces $R\in\Lfield(X)$.
We prove both this classification and the following quantitative
result independently of the analytic rigidity theorem of Part~I.

A \texttt{Lean~4} formalization of the complete equivalence in
\Cref{thm:classification} is available at
\url{https://github.com/DiegoMarquesMath/MahlerLeanII}.

For a real number $y$, we use the \textit{irrationality exponent}
\[
 \mu(y):=\sup\{t>0:0<|y-a/b|<b^{-t}
       \text{ for infinitely many }(a,b)\in\Z\times\Z_{>0}\}.
\]
Thus an irrational $y$ belongs to $\Lio$ exactly when
$\mu(y)=\infty$.

\begin{theorem}[Finite-exponent escape]\label{thm:escape}
Let $R\in\R(X)\setminus\Lfield(X)$ be nonconstant. There is a
constant $\tau_R>2$ such that every nonempty open interval
$V\subseteq D_R$ contains a Cantor set $C\subseteq V\cap\Lio$
and an integer $b_0\geq2$ for which
\begin{equation}\label{eq:escape-intro}
 \abs{R(\xi)-\frac ab}>b^{-\tau_R}
 \qquad
 (\xi\in C,\ a\in\Z,\ b\in\Z,\ b\geq b_0).
\end{equation}
In particular, every $R(\xi)$, $\xi\in C$, is irrational and has
irrationality exponent at most $\tau_R$.
One may take
\[
 \tau_R=10+2n_R(n_R-1)+2\nu_R+
        \nu_R n_R(2d_R+1),
\]
where $d_R=\max\{\deg P,\deg Q\}$ and
$n_R=\deg P+\deg Q+1$ for the normalized pair $R=P/Q$.
Here $\nu_R$ is the least obstruction order occurring among
its coefficients outside $\Lfield$, in the sense of
\Cref{lem:obstruction}.
\end{theorem}

The exponent $\tau_R$ is independent of $V$, and its dependence
on $R$ is necessary even for translations, by
\Cref{prop:no-uniform-escape}. The cutoff $b_0$ is uniform on $C$. The fractions in
\eqref{eq:escape-intro} need not be reduced, so the lower bound
excludes rational images as well as Liouville images.

Removing the rational alternative from
\Cref{thm:classification} gives the exact criterion for Maillet's
property.

\begin{corollary}[Exact preservation]\label{cor:exact}
Let $R=P/Q$ be nonconstant and normalized as in
\eqref{eq:normalized-pair}, and let
$U\subseteq D_R$ be a nonempty open interval. Then
\begin{equation}\label{eq:exact-intro}
 R(\Lio\cap U)\subseteq\Lio
 \quad\Longleftrightarrow\quad
 \begin{cases}
  P,Q\in\Lfield[X],\\
  Z_{\R}(P-tQ)\cap U\cap\Lio=\varnothing
          &\text{for every }t\in\Q.
 \end{cases}
\end{equation}
The same equivalence holds with $U$ replaced by $D_R$.
\end{corollary}

Here $Z_{\R}(A)$ denotes the real zero set of a polynomial $A$.
Coprimality ensures that a zero of $P-tQ$ is not a pole of $R$.
The coefficient condition is the principal arithmetic assertion: it
rules out every irrational non-Liouville image. The fiber condition
then excludes rational images, and is not claimed to be an effective
decision procedure. It cannot be omitted, even in higher degree.
Indeed, let $\mu\in\Lfield\setminus\Q$ and $m\geq2$. Since
$\Lfield$ is a field, $\lambda:=\mu^{-m}$ lies in
$\Lfield\setminus\Q$, yet
$\lambda\mu^m=1$. Thus $\lambda X^m\in\Lfield[X]$
fails Maillet's property at the Liouville number $\mu$.
In contrast, the coefficients constructed in
\Cref{thm:strong-nonlinear,thm:example} have no Liouville points
in their rational fibers.

Together with \cite[Corollary~1.2]{MarquesI}, this gives the
corresponding characterization for real-analytic functions.
A nonconstant real-analytic function on an open interval $U$ has
Maillet's property if and only if it agrees on $U$ with a rational
function whose normalized coefficients belong to $\Lfield$,
whose denominator has no zero in $U$, and whose rational fibers
contain no Liouville point of $U$. For entire functions, the same
criterion applies with $Q=1$. Constant functions have Maillet's
property precisely when their value belongs to $\Lio$.

Our convention requires preservation only at points of the domain.
If one additionally requires a rational function to be defined at
every Liouville number, the criterion must include
$Z_{\R}(Q)\cap\Lio=\varnothing$. The classification is structural;
it does not assert a finite decision procedure for the fiber
condition.

\subsection{Degree-one rigidity and nonlinear sharpness}

The coefficient classification distinguishes sharply between degree
one and higher degrees. In the latter case, the following result,
proved in \Cref{sec:strong}, gives much more than the existence
of a single transcendental coefficient.

\begin{theorem}[Strong Liouville coefficients]
\label{thm:strong-nonlinear}
There exist continuum many positive strong Liouville numbers
$\sigma$ such that, for every integer $m\geq2$ and every
$r\in\Q^\times$, every complex root of
\[
 X^m-\frac r\sigma
\]
is a $U_m$-number. Consequently, each such $\sigma$ satisfies
\begin{equation}\label{eq:strong-preservation}
 \sigma\xi^m\in\Lio
 \qquad(\xi\in\Lio,\ m\geq2).
\end{equation}
\end{theorem}

In particular, no algebraic coefficient restriction can be
imposed on the polynomials in the entire-function conclusion of
Part~I. The stronger fiber statement distinguishes the degree-$m$
case from the degree-one rigidity in \Cref{cor:mobius}.

\subsection{Why this field is canonical}

The approximation condition in \Cref{def:field} admits an intrinsic
interpretation through the arithmetic of $\Eset$. The classification
gives
\begin{equation}\label{eq:stabilizer-intro}
 \begin{split}
 \Lfield
 &=\{\alpha\in\R:\alpha\Lio\subseteq\Eset\}\\
 &=\{\alpha\in\R:\alpha+\Lio\subseteq\Eset\}\\
 &=\{\alpha\in\R:\Eset+\alpha=\Eset\}.
 \end{split}
\end{equation}
Moreover, $\alpha\Eset=\Eset$ for every nonzero
$\alpha\in\Lfield$. Thus the same field is determined by
multiplication of Liouville numbers, by addition to Liouville
numbers, and by the translation symmetries of $\Eset$. These
identities are proved in \Cref{sec:consequences}.

Fields consisting only of rational and Liouville numbers have been
studied before; see Marcos \cite{Marcos1999} and Kumar, Thangadurai,
and Waldschmidt \cite{KumarThangaduraiWaldschmidt2015}.
The latter authors study Liouville sets and fields through
approximation conditions along common sequences. The point
here is not the existence of such a field, but the identification of
the field that governs every nonconstant real rational function
satisfying \eqref{eq:weak-intro}. Indeed, Maillet's theorem gives
$\Q(\lambda)\subseteq\Eset$ for every $\lambda\in\Lio$, whereas
not every Liouville number belongs to $\Lfield$. Belonging to a
field contained in $\Eset$ is therefore weaker than belonging to
the coefficient field in \Cref{thm:classification}.

The field $\Lfield$ contains every strong Liouville number
(see Section~\ref{subsection:size} for the definition) and has cardinality
and transcendence degree equal to the continuum.
Nevertheless, it is meagre and has Hausdorff dimension zero.
We also compare its defining approximation property with the
fixed-gap condition considered in
\cite{BiluMarques2026, ChavesMarquesTrojovsky2021}.

The strong Liouville coefficients in
\Cref{thm:strong-nonlinear} belong to the fixed-gap class
$\Ggap$ defined in \Cref{sec:examples}. This is not the
full range of possible coefficients: we construct an explicit
positive $\alpha\in\Lfield\setminus\Q(\Ggap)$ such that
\[
 \alpha\xi^m\in\Lio
 \qquad(\xi\in\Lio,\ m\geq2).
\]

\subsection{Proof strategy}

Sufficiency follows by synchronizing finitely many coefficient
approximations and matching them to a sufficiently accurate rational approximation to a Liouville number. The denominator bound in
\eqref{eq:scale-property} controls the height of the resulting
rational approximation to the image. Unless the image is rational,
this produces approximations of arbitrarily high order.

Necessity rests on rational reconstruction. Failure of
\eqref{eq:scale-property} provides a fixed approximation order
and arbitrarily large scales at which no approximation of that
order exists within a prescribed polynomial denominator range.
At such a scale, sufficiently many rational points close to the
graph of $R$ would reconstruct a rational coefficient vector
that contradicts this obstruction.

The reconstruction uses $\deg P+\deg Q+1$ points.
A Vandermonde determinant controls its stability, while clearing
the source and target denominators separately controls its
arithmetic height. Their combination yields a two-height counting
estimate with a subquadratic error, leaving enough rational
source points from which to construct Liouville numbers whose images are irrational and have finite irrationality exponent.

The main difficulty is that the counting estimate is available
only along selected source scales: the obstruction need not persist
at every sufficiently large height. We show that arbitrarily
large admissible heights, chosen after the data of each stage
have been fixed, suffice for a nested Cantor construction.
The source approximation orders increase without bound, while
successive target ranges cover every sufficiently large
denominator.

For the strong Liouville examples, we prescribe prime valuations
of continued-fraction denominators so that the polynomials defining
nonzero rational fibers are Eisenstein. A gap principle for algebraic
approximation then places these fibers in $U_m$; the coefficient-field
criterion gives exact preservation of Liouville numbers.

The analytic rigidity theorem of Part~I enters only in the passage
to real-analytic, holomorphic, and entire functions in
\Cref{cor:analytic,cor:holomorphic,cor:entire}.

\medskip
\noindent\emph{Acknowledgments.}
The author is grateful to Yuri Bilu, Carlos Gustavo Moreira,
Johannes Schleischitz, and Jean Lelis for valuable discussions
concerning Mahler's problem and Liouville numbers.
He is especially grateful to Michel Waldschmidt, who first
introduced him to the problem at the 2008 Arizona Winter School.

\section{Approximation on every scale}
\label{sec:field}

\noindent
\textbf{Notation and conventions.}
For any finite union of real intervals, vertical bars denote
Lebesgue measure. We write $\den(r)$ for the positive denominator
of $r\in\Q$ in lowest terms. For vectors
$\mathbf x=(x_1,\ldots,x_s)\in\R^s$, we use the sup norm
$\norm{\mathbf x}_\infty:=\max_{1\leq i\leq s}\abs{x_i}$.

The symbols $O(\cdot)$, $\ll$, and $\asymp$ have their
usual meanings. Unless otherwise stated, implied constants
may depend on quantities fixed before an estimate, but not
on parameters subsequently allowed to vary.
Additional dependence is indicated by subscripts; thus
$X\ll_{A,J}Y$ means that $X\leq C(A,J)Y$ for some
constant $C(A,J)>0$. Generic constants
$C,C_1,C_2,\ldots$ may vary from one occurrence to the next.

Unless explicitly stated otherwise, scale and denominator
parameters are positive integers. When uniformity in a
parameter is asserted, both the implied constants and the
relevant scale thresholds are independent of that parameter.

\subsection{Elementary consequences and simultaneous approximation}

We first locate $\Lfield$ within $\Eset$, then show that its
approximation property is compatible with the field operations.

\begin{lemma}\label{lem:elementary-field}
One has $\Q\subseteq\Lfield\subseteq\Eset$.
\end{lemma}
\begin{proof}
The inclusion $\Q\subseteq\Lfield$ follows by using exact rational
approximations once $H$ exceeds the denominator.

Suppose now, towards a contradiction, that $\alpha\in\Lfield$ is
irrational but not Liouville. There exist constants $C,\eta>0$
such that $|\alpha-p/q|>Cq^{-\eta}$ for all $p\in\Z$ and all
integers $q\geq1$. Choose an integer $N>\eta$. By
\eqref{eq:scale-property}, for every sufficiently large integer $H$
there are integers $p_H$ and $q_H\geq1$ such that
\[
 C
 <q_H^\eta\left|\alpha-\frac{p_H}{q_H}\right|
 \leq H^{-N}q_H^{\eta-N}
 \leq H^{-N}.
\]
Letting $H\to\infty$ gives a contradiction. Thus every irrational
element of $\Lfield$ is Liouville, proving the second inclusion.
\end{proof}

To combine elements of $\Lfield$, we need approximations with a
common denominator. The next lemma provides these while retaining
polynomial control in the scale.

\begin{lemma}[Simultaneous approximation]
\label{lem:synchronization}
Let $\boldsymbol\alpha=(\alpha_1,\ldots,\alpha_s)\in\Lfield^s$.
For every integer $N\geq1$ there are integers $D\geq1$ and
$H_0\geq2$ such that every integer $H\geq H_0$ admits
$h\in\Z_{>0}$ and $\mathbf p\in\Z^s$ with
\begin{equation}\label{eq:simultaneous}
 h\leq H^D,
 \qquad
 \norm{\boldsymbol\alpha-\mathbf p/h}_\infty
       \leq(Hh)^{-N}.
\end{equation}
\end{lemma}
\begin{proof}
We argue by induction on $s$, the case $s=1$ being the definition.
Write $\boldsymbol\alpha=(\mathbf a,\beta)$ and assume the assertion
for $\mathbf a$. Fix $N$, and choose an exponent $D_a$ for
$\mathbf a$ at order $N$.

We first approximate $\beta$ at the higher order $N(D_a+1)$.
By \eqref{eq:scale-property}, there is an exponent $D_b$ such
that, for every sufficiently large $H$, we can choose $b/k$ with
\[
 k\leq H^{D_b},\qquad
 |\beta-b/k|\leq(Hk)^{-N(D_a+1)}.
\]
Apply the induction hypothesis to $\mathbf a$ at scale $Hk$,
obtaining $\mathbf c\in\Z^{s-1}$ and an integer $l\geq1$ such that
\[
 l\leq(Hk)^{D_a},\qquad
 \norm{\mathbf a-\mathbf c/l}_\infty\leq(Hkl)^{-N}.
\]
The required scale threshold is satisfied once $H$ is sufficiently
large, since $k\geq1$.

Set $h=kl$. The approximation to $\mathbf a$ already has the
desired accuracy, and $Hh=Hkl\leq(Hk)^{D_a+1}$ gives
$|\beta-b/k|\leq(Hh)^{-N}$. Moreover,
\[
 h\leq H^{D_b+D_a(D_b+1)}.
\]
Thus $\mathbf p=(k\mathbf c,lb)$ and this choice of $h$ satisfy
\eqref{eq:simultaneous}, completing the induction.
\end{proof}

The common denominator allows us to control the approximation
errors in sums and products. Inversion can be handled directly.

\begin{proposition}[Field structure]\label{prop:field}
The set $\Lfield$ is a subfield of $\R$. In particular, for every
integer $s\geq1$ and every $F\in\Q(X_1,\ldots,X_s)$, one has
$F(\alpha_1,\ldots,\alpha_s)\in\Lfield$ whenever
$\alpha_1,\ldots,\alpha_s\in\Lfield$ and $F$ is defined at
$(\alpha_1,\ldots,\alpha_s)$.
\end{proposition}
\begin{proof}
By \Cref{lem:elementary-field}, $\Q\subseteq\Lfield$. It therefore
suffices to prove closure under subtraction, multiplication, and
inversion.

Let $\alpha,\beta\in\Lfield$ and fix an integer $N\geq1$.
By \Cref{lem:synchronization}, there is an integer $D\geq1$
such that, for every sufficiently large integer $H$, we can find
$a,b\in\Z$ and $h\in\Z_{>0}$ with $h\leq H^D$ and
\[
 \max\left\{
 \abs{\alpha-\frac ah},
 \abs{\beta-\frac bh}
 \right\}\leq(Hh)^{-(2N+1)}.
\]
In particular, $|b/h|\leq|\beta|+1$ for sufficiently large $H$.
The triangle inequality then gives
\[
 \begin{aligned}
 \abs{\alpha-\beta-\frac{a-b}{h}}
 &\leq 2(Hh)^{-(2N+1)}
 \leq(Hh)^{-N},\\
 \abs{\alpha\beta-\frac{ab}{h^2}}
 &\leq (|\alpha|+|\beta|+1)(Hh)^{-(2N+1)}
 \leq(Hh^2)^{-N},
 \end{aligned}
\]
where the last inequalities hold after increasing the threshold
for $H$. The denominators are bounded by $H^D$ and $H^{2D}$,
respectively. Since $N$ was arbitrary, this proves that
$\alpha-\beta,\alpha\beta\in\Lfield$.

Now let $\alpha\in\Lfield\setminus\{0\}$ and fix an integer
$N\geq1$. Applying \eqref{eq:scale-property} at order $N+1$,
we obtain an integer $D\geq1$ such that every sufficiently large
$H$ admits $a\in\Z$ and $h\in\Z_{>0}$ satisfying
$h\leq H^D$ and
$|\alpha-a/h|\leq(Hh)^{-(N+1)}$.
For sufficiently large $H$, we have
$|a/h|\geq|\alpha|/2$ and $|a|\leq Ah$, where
$A:=\max\{1,2|\alpha|\}$. In particular, $a\ne0$, and
\[
 \abs{\alpha^{-1}-\frac ha}
 =\frac{|\alpha-a/h|}{|\alpha|\,|a/h|}
 \leq\frac{2}{|\alpha|^2}(Hh)^{-(N+1)}
 \leq(H|a|)^{-N}.
\]
The last inequality holds once $H\geq2A^N/|\alpha|^2$.
Written with positive denominator, $h/a$ has denominator
$|a|\leq AH^D\leq H^{D+1}$ for sufficiently large $H$.
Thus $\alpha^{-1}$ satisfies \eqref{eq:scale-property} at
every order, and hence belongs to $\Lfield$.

This proves that $\Lfield$ is a field. The assertion about
rational functions follows by evaluating their numerators
and nonzero denominators in this field.
\end{proof}

\begin{remark}\label{rem:coefficient-normalization}
For the normalized pair in \eqref{eq:normalized-pair},
$R\in\Lfield(X)$ if and only if $P,Q\in\Lfield[X]$.
Indeed, a coprime pair over $\Lfield$ remains coprime over $\R$
by B\'ezout's identity, so its monic normalization over $\Lfield$
is also the unique normalized real pair.
Any other reduced real pair $(A,B)$ representing $R$ differs
from $(P,Q)$ by a common scalar. The equivalent condition is
therefore $A/c,B/c\in\Lfield[X]$, where $c$ is any nonzero
coefficient of $(A,B)$.
\end{remark}

\subsection{Size and nontriviality}\label{subsection:size}

A real irrational $\sigma$ is \emph{strong Liouville} if the
denominators $q_j$ of all its continued-fraction convergents satisfy
\begin{equation}\label{eq:strong}
 \frac{\log q_{j+1}}{\log q_j}\longrightarrow\infty \qquad (j\to \infty).
\end{equation}
This is equivalent to the usual definition by approximation orders
tending to infinity along all convergents; see
\cite{Bugeaud2004,Petruska1992}.

\begin{proposition}\label{prop:strong}
The field $\Lfield$ contains every strong Liouville number.
In particular, $|\Lfield|=\operatorname{trdeg}_{\Q}\Lfield=2^{\aleph_0}$.
\end{proposition}
\begin{proof}
Let $\sigma$ be a strong Liouville number, and let $(p_j/q_j)$
be the convergents of its continued fraction. Fix an integer $N\geq1$.
By \eqref{eq:strong}, we may choose $j_0$ such
that $q_{j+1}\geq q_j^{2N}$ for every $j\geq j_0$.
For each sufficiently large integer $H$, choose $j\geq j_0$
such that $q_j\leq H<q_{j+1}$.

If $q_{j+1}\geq H^{2N}$, the convergent $p_j/q_j$ satisfies
\[
 \abs{\sigma-\frac{p_j}{q_j}}
 <\frac{1}{q_jq_{j+1}}
 \leq H^{-2N}
 \leq(Hq_j)^{-N},
\]
and its denominator is at most $H$.
If $q_{j+1}<H^{2N}$, we instead use $p_{j+1}/q_{j+1}$.
Since $q_{j+2}\geq q_{j+1}^{2N}$ and $H<q_{j+1}$, we have
\[
 \abs{\sigma-\frac{p_{j+1}}{q_{j+1}}}
 <\frac{1}{q_{j+1}q_{j+2}}
 \leq q_{j+1}^{-(2N+1)}
 <(Hq_{j+1})^{-N}.
\]
In either case, we obtain an approximation satisfying
\eqref{eq:scale-property} with denominator at most $H^{2N}$.
Thus $D=2N$ is admissible for every $N$, proving that
$\sigma\in\Lfield$.

At each stage of a continued-fraction expansion, choose between
two distinct positive partial quotients ensuring
$q_{j+1}>q_j^{j+1}$. The resulting $2^{\aleph_0}$ expansions
represent distinct strong Liouville numbers, so the inclusion
just proved gives $|\Lfield|=2^{\aleph_0}$.

If $T$ is a transcendence basis of $\Lfield/\Q$, then
$\Lfield$ is algebraic over $\Q(T)$, and hence
\[
 2^{\aleph_0}=|\Lfield|=|\Q(T)|
             =\max\{\aleph_0,|T|\}.
\]
It follows that $|T|=2^{\aleph_0}$.
\end{proof}

We now turn to the topological and metric properties of $\Lfield$.

\begin{proposition}\label{prop:size}
The field $\Lfield$ is a dense Borel subset of $\R$, is meagre,
and has Hausdorff dimension zero. Moreover,
$\Lio\setminus\Lfield$ is residual in $\R$; in particular,
$\Lfield\subsetneq\Eset$.
\end{proposition}
\begin{proof}
Density follows from $\Q\subseteq\Lfield$. The
Jarn\'ik--Besicovitch theorem gives $\dim_H\Lio=0$
\cite[Theorem~5.2]{Bugeaud2004}. Since $\Q$ is countable and
$\Lfield\subseteq\Q\cup\Lio$, it follows that
$\dim_H\Lfield=0$.

For fixed integers $N,D,H$, the real numbers satisfying
\eqref{eq:scale-property} for some $p\in\Z$ and $1\leq q\leq H^D$
form a countable union of closed intervals. The remaining
quantifiers in the definition of $\Lfield$ involve only countable
unions and intersections, so $\Lfield$ is Borel.

Suppose that $\Lfield$ were nonmeagre. Its Baire property would
make it comeagre in some nonempty open interval $I$. For every
sufficiently small $t$, the sets $\Lfield$ and $\Lfield+t$ would
then meet in $I\cap(I+t)$, giving
$t\in\Lfield-\Lfield=\Lfield$. Thus $\Lfield$ would contain a
neighbourhood of zero and, being an additive subgroup, would
equal $\R$, contrary to $\dim_H\Lfield=0$.

Finally, $\Lio$ is a dense $G_\delta$ subset of $\R$. Removing
the meagre set $\Lfield$ leaves a residual set, proving the
remaining assertions.
\end{proof}

\section{Values of rational functions at Liouville numbers}
\label{sec:sufficiency}

We prove the sufficiency direction of \Cref{thm:classification}.
Given a rational approximation $p/b$ to a Liouville number $\xi$,
we use $b$ as the scale for approximating the coefficients of $R$.

\begin{proposition}\label{prop:evaluation}
Every $R\in\Lfield(X)$ satisfies $R(\Lio\cap D_R)\subseteq\Eset$.
\end{proposition}

\begin{proof}
Write $R=P/Q$ in reduced form over $\Lfield$, and choose an integer
$d\geq1$ such that $P$ and $Q$ have degree at most $d$. Fix
$\xi\in\Lio\cap D_R$ and put $z=R(\xi)$. We may assume that
$z$ is irrational.

Fix an integer $k\geq1$, and let $D$ be the exponent given by
\Cref{lem:synchronization} for the coefficients of $P$ and $Q$
at order $d(k+1)$. Since $\xi$ is Liouville, there is a sequence
of reduced fractions $p/b$ with $b\to\infty$ and
$|\xi-p/b|<b^{-(D+d)(k+1)}$. For each sufficiently large $b$ in
this sequence, the lemma at scale $b$ gives polynomials
$\widetilde P,\widetilde Q\in h^{-1}\Z[X]$ of degree at most $d$,
with $1\leq h\leq b^D$, whose coefficients differ from those of
$P,Q$, respectively, by at most $(bh)^{-d(k+1)}$.

Set $T_b=hb^d$. Since $T_b\leq b^{D+d}$ and $T_b\leq(bh)^d$,
both $|\xi-p/b|$ and the coefficient errors are at most
$T_b^{-(k+1)}$. The quotient $A(x)/B(x)$ depends smoothly on
$x$ and on the coefficients of $A$ and $B$ near $(P,Q,\xi)$,
because $Q(\xi)\ne0$. Thus, for all sufficiently large $b$, the
rational number $z_b=\widetilde P(p/b)/\widetilde Q(p/b)$ is
well defined and satisfies
\begin{equation}\label{eq:evaluation-error}
 |z-z_b|\ll T_b^{-(k+1)}.
\end{equation}

The factor $T_b=hb^d$ clears the denominators of both
$\widetilde P(p/b)$ and $\widetilde Q(p/b)$. Since
$\widetilde Q(p/b)\to Q(\xi)\ne0$, we obtain
\begin{equation}\label{eq:evaluation-height}
 \den(z_b)\leq T_b\abs{\widetilde Q(p/b)}\ll T_b.
\end{equation}
Combining \eqref{eq:evaluation-error} and \eqref{eq:evaluation-height}
gives $\den(z_b)^k|z-z_b|\ll T_b^{-1}\to0$. As the rationals
$z_b$ converge to the irrational number $z$, their reduced
denominators tend to infinity. Hence $0<|z-z_b|<\den(z_b)^{-k}$
for arbitrarily large denominators. Since $k$ was arbitrary,
$z\in\Lio$.
\end{proof}

\begin{remark}\label{rem:extension-field}
For each $\xi\in\Lio$, the values $R(\xi)$ with $R\in\Lfield(X)$
and $\xi\in D_R$ form precisely the field $\Lfield(\xi)$.
Thus $\Lfield(\xi)\subseteq\Eset$. When $\xi\notin\Lfield$,
this gives a proper extension of $\Lfield$ inside $\Eset$.
For such a $\xi$, \Cref{cor:intrinsic} will show that
$\xi+\eta\notin\Eset$ for some $\eta\in\Lio$.
\end{remark}

\section{Rational reconstruction}
\label{sec:reconstruction}

We establish stable reconstruction for rational functions of
arbitrary bidegree, including the three-point M\"obius case.
The separation of the abscissae controls the error in
the reconstructed coefficients. When the points are rational,
the denominators of their coordinates also give a bound for a
common denominator of these coefficients.

Fix a nonconstant $R=P/Q\in\R(X)$ in the normalization
\eqref{eq:normalized-pair}, and set
\begin{equation}\label{eq:degree-parameters}
 \begin{gathered}
 r:=\deg P,\qquad s:=\deg Q,\qquad d:=\max\{r,s\},\\
 n:=r+s+1,\qquad k:=n(n-1).
 \end{gathered}
\end{equation}
Thus $d\geq1$ and $n\geq2$. Let $\boldsymbol\theta\in\R^n$
be the column vector consisting of the coefficients of $P$ in
degrees $0,\ldots,r$, followed by those of $Q$ in degrees
$0,\ldots,s-1$. Zero coefficients are included. The leading
coefficient of $Q$ is fixed at $1$.

\begin{lemma}[Stable rational interpolation]
\label{lem:reconstruction}
Let $I_*\subseteq D_R$ be a nondegenerate compact interval.
There are constants $c_0,C_0>0$, depending only on $R$ and
$I_*$, with the following property. For distinct
$x_1,\ldots,x_n\in I_*$ and arbitrary $y_1,\ldots,y_n\in\R$, set
\[
 \delta:=\max_{1\leq i\leq n}|y_i-R(x_i)|,
 \qquad
 V(\mathbf x):=\prod_{1\leq i<j\leq n}|x_i-x_j|.
\]
If $\delta\leq c_0V(\mathbf x)$, there exist polynomials
$\widehat P,\widehat Q\in\R[X]$, with $\deg\widehat P\leq r$
and $\widehat Q$ monic of degree $s$, such that $\widehat Q$
has no zero on $I_*$ and
\[
 \frac{\widehat P(x_i)}{\widehat Q(x_i)}=y_i
 \qquad(1\leq i\leq n).
\]
Their free coefficient vector $\widehat{\boldsymbol\theta}$
satisfies
\begin{equation}\label{eq:reconstruction-error}
 \norm{\widehat{\boldsymbol\theta}-\boldsymbol\theta}_\infty
 \leq C_0\frac{\delta}{V(\mathbf x)}.
\end{equation}

Suppose, in addition, that $x_i=p_i/q_i$ and $y_i=a_i/b_i$,
with integer numerators and positive integer denominators
satisfying $q_i<2S$ and $b_i<2B$, where $S,B\geq1$.
Then $\widehat{\boldsymbol\theta}=\mathbf c/h$ for some
$\mathbf c\in\Z^n$ and $h\in\Z_{>0}$ with
\begin{equation}\label{eq:reconstruction-height}
 h\leq C_0 S^{nd}B^n.
\end{equation}
\end{lemma}

\begin{proof}
Since the leading coefficient of $\widehat Q$ is fixed, the
equations $\widehat P(x_i)=y_i\widehat Q(x_i)$ form a square
linear system $M\widehat{\boldsymbol\theta}=\mathbf t$.
The $i$th row of $M$ and the corresponding entry of $\mathbf t$
are
\[
 \bigl((x_i^j)_{0\leq j\leq r},
       (-y_i x_i^j)_{0\leq j<s}\bigr),
 \qquad t_i=y_i x_i^s.
\]
When $s=0$, the second block is empty. Let $M_0$ denote the
matrix obtained by replacing each $y_i$ with $R(x_i)$.

Consider the linear map
\begin{equation}\label{eq:polynomial-map}
 \R[X]_{\leq r}\times\R[X]_{<s}\longrightarrow\R[X]_{\leq n-1},
 \qquad (U,W)\longmapsto QU-PW,
\end{equation}
where $\R[X]_{<s}=\{0\}$ if $s=0$. If $QU=PW$, coprimality of
$P$ and $Q$ gives $Q\mid W$. Since $\deg W<s=\deg Q$, we have
$W=0$ and hence $U=0$. Both spaces have dimension $n$, so
\eqref{eq:polynomial-map} is an isomorphism. In particular, the
polynomials $X^jQ$ for $0\leq j\leq r$ and $-X^jP$ for
$0\leq j<s$ form a basis of $\R[X]_{\leq n-1}$.

Multiplying the $i$th row of $M_0$ by $Q(x_i)$ gives the values
of these basis polynomials at $x_i$. The Vandermonde determinant
therefore yields
\begin{equation}\label{eq:vandermonde-minor}
 \det M_0
 =c_R\,\frac{\prod_{i<j}(x_j-x_i)}{\prod_i Q(x_i)}
\end{equation}
for some nonzero constant $c_R$. Since $Q$ is bounded on $I_*$,
it follows that $|\det M_0|\gg V(\mathbf x)$. Here and below,
the implied constants depend only on $R$ and $I_*$.

The product $V(\mathbf x)$ is bounded on $I_*^n$, so we may
choose $c_0$ small enough that $\delta\leq1$. The entries of
$M$ and $M_0$ are then uniformly bounded and differ by
$O(\delta)$. Multilinearity of the determinant gives
$|\det M-\det M_0|\ll\delta$. Thus, for $c_0$ sufficiently small,
\[
 |\det M|\geq\tfrac12|\det M_0|\gg V(\mathbf x).
\]
In particular, $M\widehat{\boldsymbol\theta}=\mathbf t$ has
a unique solution.

For the original coefficient vector $\boldsymbol\theta$, we have
\[
 (\mathbf t-M\boldsymbol\theta)_i
 =Q(x_i)\bigl(y_i-R(x_i)\bigr).
\]
Each entry is $O(\delta)$. Applying Cramer's rule to
$M(\widehat{\boldsymbol\theta}-\boldsymbol\theta)
=\mathbf t-M\boldsymbol\theta$, and using the uniform bound on
the cofactors of $M$, gives \eqref{eq:reconstruction-error}.
Consequently,
\[
 \sup_{x\in I_*}|\widehat Q(x)-Q(x)|
 \ll \frac{\delta}{V(\mathbf x)}\leq c_0.
\]
As $\min_{x\in I_*}|Q(x)|>0$, a further decrease of $c_0$
ensures that $\widehat Q$ has no zero on $I_*$. The defining
equations now give the asserted interpolation identities.

Finally, suppose that $x_i=p_i/q_i$ and $y_i=a_i/b_i$ as in the
statement. Multiplying the $i$th equation by $q_i^d b_i$ clears
all denominators, including that of $t_i$. The resulting matrix
$M_{\Z}$ and right-hand side have integer entries of size
$O(S^dB)$, since the $x_i$ and $y_i$ are uniformly bounded.
Thus
\[
 0<|\det M_{\Z}|\ll S^{nd}B^n.
\]
Cramer's rule gives the common denominator $h=|\det M_{\Z}|$,
which proves \eqref{eq:reconstruction-height}.
\end{proof}

\section{Two-height counting at obstructed scales}
\label{sec:counting}

If $R\notin\Lfield(X)$, one of its normalized coefficients fails
the approximation condition at arbitrarily large scales. At these
scales, \Cref{lem:reconstruction} bounds the number of rational
points near the graph of $R$ whose two denominators lie in
prescribed ranges.

For $\alpha\in\R$ and positive integers $\nu,D$, define
\begin{equation}\label{eq:obstructed-set}
\begin{split}
 \mathcal O_D(\alpha,\nu):=\bigl\{S\in\Z_{\geq2}:\;&
 \abs{\alpha-a/h}>(Sh)^{-\nu}\\
 &\text{for all }a\in\Z,\quad h=1,\ldots,S^D\bigr\}.
\end{split}
\end{equation}

\begin{lemma}[Obstructed scales]\label{lem:obstruction}
A real number $\alpha$ lies outside $\Lfield$ if and only if
there is an integer $\nu\geq1$ such that
$\mathcal O_D(\alpha,\nu)$ is unbounded for every integer $D\geq1$.
\end{lemma}
\begin{proof}
By \Cref{def:field}, $\alpha\notin\Lfield$ precisely when there
is an order $\nu\geq1$ for which the following holds. Given any
integers $D\geq1$ and $S_0\geq2$, there is some $S\geq S_0$
at which no fraction $a/h$ with $1\leq h\leq S^D$ satisfies
$|\alpha-a/h|\leq(Sh)^{-\nu}$. Equivalently,
$\mathcal O_D(\alpha,\nu)$ is unbounded for every $D$.
\end{proof}

\subsection{Counting reduced fractions in intervals}

For each positive integer $S$, write
\[
 \mathcal F_S:=
 \{p/q:p\in\Z,\ q\in\Z_{>0},\ \gcd(p,q)=1,\ S\leq q<2S\}.
\]

\begin{lemma}\label{lem:farey}
If $E\subseteq\R$ is a union of at most $M$ intervals, then
\begin{equation}\label{eq:farey-upper}
 \#(E\cap\mathcal F_S)\leq 4S^2|E|+M.
\end{equation}
For every nondegenerate compact interval $J$,
\begin{equation}\label{eq:farey-lower}
 \#(J\cap\mathcal F_S)\geq \cF |J|S^2,
 \qquad \cF=\tfrac14,
\end{equation}
for all sufficiently large $S$, with the threshold depending on $J$.
\end{lemma}
\begin{proof}
Distinct fractions in $\mathcal F_S$ differ by more than
$(4S^2)^{-1}$. An interval of length $\ell$ therefore contains
at most $4S^2\ell+1$ such fractions. Summing over the components
of $E$ gives \eqref{eq:farey-upper}.

For $J=[a,b]$, M\"obius inversion gives
\[
 \#\{p\in\Z:aq\leq p\leq bq,\ \gcd(p,q)=1\}
 =|J|\varphi(q)+O(\tau_0(q)),
\]
where $\tau_0(q)$ is the number of positive divisors of $q$.
Summing over $S\leq q<2S$ and applying the classical summatory
estimates for $\varphi$ and $\tau_0$ \cite{Apostol1976}, we obtain
\[
 \#(J\cap\mathcal F_S)
 =\frac{9}{\pi^2}|J|S^2+O\bigl((1+|J|)S\log S\bigr).
\]
Since $9/\pi^2>1/4$, this proves \eqref{eq:farey-lower}.
\end{proof}

\subsection{The counting estimate}

Assume henceforth that $R\notin\Lfield(X)$ and retain the notation
of \eqref{eq:degree-parameters}. By
\Cref{rem:coefficient-normalization}, some coordinate of
$\boldsymbol\theta$ lies outside $\Lfield$. Choose such a
coordinate $\alpha$ whose obstruction order is $\nu_R$,
as in \Cref{thm:escape,lem:obstruction}, and set
\begin{equation}\label{eq:tau-choice}
 \tau:=10+2k+2\nu+\nu n(2d+1).
\end{equation}
In particular, $\tau>3$, and its choice depends only on $R$.

Fix a nondegenerate compact interval $I_*\subseteq D_R$ on which
\begin{equation}\label{eq:derivative-bounds}
 0<m_*\leq |R'(x)|\leq M_*
 \qquad(x\in I_*).
\end{equation}
Then $R$ is strictly monotone on $I_*$, and the mean value
theorem gives $|x-x'|\leq m_*^{-1}|R(x)-R(x')|$ there.

For integers $A\geq3$ and $S\geq2$, put
\begin{equation}\label{eq:target-height}
 T(S,A):=S^{A/(\tau-3)}.
\end{equation}
Given $J\subseteq I_*$ and an integer $H\geq2$, define
\begin{equation}\label{eq:resonant-set}
\begin{split}
 \mathcal D(J;S,A,H):=\bigl\{x\in J\cap\mathcal F_S:\;&
 |R(x)-a/b|\leq 2b^{-\tau}+4M_*S^{-A}\\
 &\text{for some }a\in\Z,\ b\in\Z_{>0},\
 H\leq b<T(S,A)\bigr\}.
\end{split}
\end{equation}
The fractions $a/b$ need not be reduced. The term $4M_*S^{-A}$
allows for the variation of $R$ over intervals of length $S^{-A}$
in the Cantor construction. At each stage of that construction,
$A$ is fixed before $S$ is chosen, so the available scales may
depend on $A$.

\begin{proposition}[Counting at obstructed scales]
\label{prop:counting}
There is a constant $\Clow>0$, depending only on $R,I_*,\tau$,
with the following property. For every integer $A\geq3$, there
are an unbounded set $\mathcal S_A\subseteq\Z_{\geq2}$ and a
constant $C_A>0$ such that
\begin{equation}\label{eq:counting}
 \#\mathcal D(J;S,A,H)
 \leq \Clow S^2H^{2-\tau}+C_AS
\end{equation}
for every $S\in\mathcal S_A$, every $J\subseteq I_*$, and every
integer $H\geq2$.
\end{proposition}
\begin{proof}
Fix $A\geq3$, and set
\begin{equation}\label{eq:denominator-range}
 D_A:=\left\lceil nd+\frac{nA}{\tau-3}\right\rceil+1.
\end{equation}
We prove the estimate for all sufficiently large
$S\in\mathcal O_{D_A}(\alpha,\nu)$. Since
$\mathcal D(J;S,A,H)=J\cap\mathcal D(I_*;S,A,H)$, it suffices
to count on $I_*$.

If $H\geq T(S,A)$, there is nothing to prove. Otherwise, divide
the range $H\leq b<T(S,A)$ into blocks
\[
 B\leq b<\min\{2B,T(S,A)\},
 \qquad B=2^jH<T(S,A),\quad j\geq0.
\]
Let $\mathcal D_B$ consist of those $x\in I_*\cap\mathcal F_S$
for which the inequality in \eqref{eq:resonant-set} holds for
some $a/b$ in this block. It suffices to sum $\#\mathcal D_B$
over the left endpoints $B$.

\smallskip\noindent\emph{Blocks with $B<\sqrt S$.}
Since $R$ is bounded on $I_*$, every fraction $a/b$ that occurs
satisfies $|a|\ll b$. Thus there are $O(B^2)$ possible pairs
$(a,b)$ in the block. For each such pair, strict monotonicity
and \eqref{eq:derivative-bounds} show that
\[
 \{x\in I_*:|R(x)-a/b|\leq 2b^{-\tau}+4M_*S^{-A}\}
\]
is either empty or an interval of length
$O(b^{-\tau}+S^{-A})$, possibly a singleton. The union of these
intervals has at most $O(B^2)$ components and total length
$O(B^{2-\tau}+B^2S^{-A})$. It follows from
\eqref{eq:farey-upper} that
\[
 \#\mathcal D_B
 \ll S^2B^{2-\tau}+S^{2-A}B^2+B^2.
\]
All implied constants here depend only on $R,I_*,\tau$.
Since $\tau>2$, the sum of $B^{2-\tau}$ over the dyadic blocks
is $O_\tau(H^{2-\tau})$. Also,
$\sum_{B<\sqrt S}B^2\ll S$ and $A\geq3$. Hence
\begin{equation}\label{eq:small-targets}
 \sum_{B<\sqrt S}\#\mathcal D_B
 \leq \Clow S^2H^{2-\tau}+O(S),
\end{equation}
with $\Clow$ independent of $A$ and $H$.

\smallskip\noindent\emph{Blocks with $B\geq\sqrt S$.}
Write $B=S^u$, so that
\begin{equation}\label{eq:u-range}
 \frac12\leq u<\frac{A}{\tau-3}.
\end{equation}
Suppose that $\mathcal D_B$ contains $n$ distinct points
$x_i=p_i/q_i$. For each $i$, choose $y_i=a_i/b_i$ satisfying
\eqref{eq:resonant-set} with $B\leq b_i<2B$. Then
\[
 \delta:=\max_i|y_i-R(x_i)|
 \leq (2+4M_*)S^{-t},
 \qquad t:=\min\{\tau u,A\}.
\]
Both $\tau u$ and $A$ are at least $(\tau-3)u$, so
\begin{equation}\label{eq:t-lower}
 t\geq(\tau-3)u.
\end{equation}
The separation of the fractions $x_i\in\mathcal F_S$ gives
\begin{equation}\label{eq:vandermonde-separation}
 V(\mathbf x)>(4S^2)^{-n(n-1)/2}=4^{-k/2}S^{-k}.
\end{equation}
Since $(\tau-3)/2>k$, we obtain
$\delta/V(\mathbf x)\ll S^{k-(\tau-3)/2}\to0$, uniformly in
the block and the choice of the $x_i$. For sufficiently large
$S$, \Cref{lem:reconstruction} therefore gives
$\mathbf c\in\Z^n$ and $h\in\Z_{>0}$ such that
\begin{equation}\label{eq:reconstructed-scale-bounds}
 \norm{\boldsymbol\theta-\mathbf c/h}_\infty\leq C_2S^{k-t},
 \qquad h\leq C_3S^{n(d+u)}.
\end{equation}
By \eqref{eq:u-range} and \eqref{eq:denominator-range}, the latter
bound is at most $S^{D_A}$ for all sufficiently large $S$,
uniformly in $u$.

Let $c$ be the coordinate of $\mathbf c$ corresponding to
$\alpha$. Since $S\in\mathcal O_{D_A}(\alpha,\nu)$ and
$h\leq S^{D_A}$, the obstruction and
\eqref{eq:reconstructed-scale-bounds} give
\begin{equation}\label{eq:contradiction-exponent}
 1<(Sh)^\nu\abs{\alpha-c/h}
 \leq C_2C_3^\nu S^{k-t+\nu+\nu n(d+u)}.
\end{equation}
On the other hand, \eqref{eq:t-lower} and \eqref{eq:tau-choice}
give
\begin{align}
 t-k-\nu-\nu n(d+u)
 &\geq (\tau-3-\nu n)u-k-\nu-\nu nd \notag\\
 &=\frac72+(\tau-3-\nu n)\bigl(u-\tfrac12\bigr)
 \geq\frac72,
 \label{eq:exponent-gap}
\end{align}
since $\tau-3-\nu n>0$ and $u\geq1/2$. The right-hand side of
\eqref{eq:contradiction-exponent} is therefore at most
$C_2C_3^\nu S^{-7/2}$, a contradiction for sufficiently large
$S$. Thus every block with $B\geq\sqrt S$ contains at most
$n-1$ points of $\mathcal F_S$ satisfying the required inequality.

There are at most $1+A\log_2 S/(\tau-3)$ blocks in total, so the
blocks with $B\geq\sqrt S$ contribute $O_A(\log S)$. Together
with \eqref{eq:small-targets}, this proves \eqref{eq:counting}.
All lower thresholds on $S$ are independent of $H$ and $J$.
Removing the scales below these thresholds from
$\mathcal O_{D_A}(\alpha,\nu)$ leaves an unbounded set
$\mathcal S_A$ by \Cref{lem:obstruction}.
\end{proof}

\section{A Cantor construction at sparse scales}
\label{sec:fusion}

We use the counting estimate to construct Liouville numbers whose
images have bounded irrationality exponent. At each generation,
the order of approximation to the points of the Cantor set
increases, while the denominators excluded for their images form
the next in a sequence of consecutive ranges.

\begin{proposition}[Escape at sparse scales]\label{prop:fusion}
Let $f$ be a real-valued $C^1$ function on a neighborhood of a
nondegenerate compact interval $I_*$, with $|f'|\leq M_*$ on
$I_*$. Fix $\tau>3$. Define $T(S,A)$ by \eqref{eq:target-height}
and $\mathcal D(J;S,A,H)$ by \eqref{eq:resonant-set}, with $f$
in place of $R$.

Suppose there is a constant $\Clow>0$ such that, for every
integer $A\geq3$, there are an unbounded set
$\mathcal S_A\subseteq\Z_{\geq2}$ and a constant $C_A>0$ for
which \eqref{eq:counting} holds for every nondegenerate compact
interval $J\subseteq I_*$, every integer $H\geq2$, and every
$S\in\mathcal S_A$. Then there exist a Cantor set
$C\subseteq\operatorname{int}(I_*)\cap\Lio$ and an integer
$b_0\geq2$ such that
\begin{equation}\label{eq:fusion-conclusion}
 \abs{f(\xi)-\frac ab}>b^{-\tau}
 \qquad
 (\xi\in C,\ a\in\Z,\ b\in\Z_{\geq b_0}).
\end{equation}
\end{proposition}

\begin{proof}
Fix an integer $A\geq3$, a nondegenerate compact interval
$J\subseteq I_*$, and an integer $H\geq2$ such that
\begin{equation}\label{eq:tail-condition}
 \Clow H^{2-\tau}\leq\frac{\cF}{4}|J|.
\end{equation}
By \eqref{eq:farey-lower} and \eqref{eq:counting}, for all
sufficiently large $S\in\mathcal S_A$ we have
\[
 \#(J\cap\mathcal F_S)-\#\mathcal D(J;S,A,H)
 \geq \frac{3\cF}{4}|J|S^2-C_AS>0.
\]
Hence there is a reduced fraction $r=p/q\in J$, with
$S\leq q<2S$, such that
\begin{equation}\label{eq:safe-center}
 \abs{f(r)-\frac ab}>2b^{-\tau}+4M_*S^{-A}
 \qquad(a\in\Z,\ H\leq b<T(S,A)).
\end{equation}
Here and throughout the proof, $b$ ranges over positive integers.

For a compact interval $I$ of positive length, let $J(I)$ denote
its closed middle third. We construct finite families
$\mathcal I_j$ of pairwise disjoint compact intervals and integers
$T_j\geq2$ such that
\begin{equation}\label{eq:tree-tail}
 \Clow T_j^{2-\tau}\leq\frac{\cF}{4}|J(I)|
 \qquad(I\in\mathcal I_j).
\end{equation}
Begin with a single nondegenerate compact interval
$I_0\subseteq\operatorname{int}(I_*)$. Set
$\mathcal I_0=\{I_0\}$ and choose $T_0$ sufficiently large for
\eqref{eq:tree-tail} to hold.

Suppose that $\mathcal I_j$ and $T_j$ have been chosen, and put
$A_j=j+3$. For every sufficiently large $S\in\mathcal S_{A_j}$,
each $J(I)$ with $I\in\mathcal I_j$ contains a fraction satisfying
\eqref{eq:safe-center} with $A=A_j$ and $H=T_j$. A single
threshold for $S$ suffices because $\mathcal I_j$ is finite.

For this fixed $A_j$, write
\[
 \rho(S)=(2S)^{-A_j},
 \qquad T'(S)=\lfloor T(S,A_j)\rfloor.
\]
Since $T'(S)\sim S^{A_j/(\tau-3)}$ as $S\to\infty$, we have
\[
 \frac{T'(S)^{2-\tau}}{\rho(S)}
 \sim 2^{A_j}S^{-A_j/(\tau-3)}\longrightarrow0.
\]
We may therefore choose $S_j\in\mathcal S_{A_j}$ above the common
threshold so that, with $\rho_j=\rho(S_j)$ and $T_{j+1}=T'(S_j)$,
\begin{align}
 S_j&>2S_{j-1}\qquad(j\geq1),\label{eq:tree-height}\\
 \rho_j&<\tfrac13\min_{I\in\mathcal I_j}|I|,
       \label{eq:tree-containment}\\
 T_{j+1}&\geq2T_j,\label{eq:tree-cutoff}\\
 \Clow T_{j+1}^{2-\tau}&\leq\frac{\cF}{48}\rho_j.
       \label{eq:tree-next-tail}
\end{align}

For each $I\in\mathcal I_j$, choose a fraction
$r_I=p_I/q_I\in J(I)\cap\mathcal F_{S_j}$ satisfying
\eqref{eq:safe-center}. Take as its two child intervals
\begin{equation}\label{eq:children}
 I^-=[r_I-\rho_j/2,r_I-\rho_j/4],
 \qquad
 I^+=[r_I+\rho_j/4,r_I+\rho_j/2].
\end{equation}
They are disjoint and lie in the interior of $I$ by
\eqref{eq:tree-containment}. Each has length $\rho_j/4$, so
its middle third has length $\rho_j/12$. Thus
\eqref{eq:tree-next-tail} gives \eqref{eq:tree-tail} at
generation $j+1$. Taking all such children as $\mathcal I_{j+1}$
completes the induction.

Two estimates hold throughout each child. First, since
$q_I<2S_j$, every $x\in I^-\cup I^+$ satisfies
\begin{equation}\label{eq:tree-source}
 0<|x-r_I|<\rho_j<q_I^{-A_j}.
\end{equation}
Second, if $T_j\leq b<T_{j+1}$, then $b<T(S_j,A_j)$, so
\eqref{eq:safe-center} applies. The derivative bound gives
$|f(x)-f(r_I)|\leq M_*\rho_j$, and hence
\begin{equation}\label{eq:tree-target}
 \abs{f(x)-\frac ab}
 >2b^{-\tau}+4M_*S_j^{-A_j}-M_*\rho_j
 >b^{-\tau}
 \qquad(a\in\Z).
\end{equation}

Set
\[
 C=\bigcap_{j\geq0}\ \bigcup_{I\in\mathcal I_j}I.
\]
The intervals at generation $j+1$ have length $\rho_j/4\to0$.
Each infinite sequence of binary choices therefore determines a
unique point of $C$, and the disjointness of the children makes
this correspondence one-to-one. It is continuous by the same
diameter bound, so $C$ is homeomorphic to $\{0,1\}^{\N}$.

Fix $\xi\in C$. Along its branch, \eqref{eq:tree-source} gives
rational approximations with denominators $q_I\geq S_j\to\infty$
and orders $A_j\to\infty$. Thus $\xi\in\Lio$. By
\eqref{eq:tree-cutoff}, the ranges $[T_j,T_{j+1})$ together
cover every integer $b\geq T_0$. For each such $b$, applying
\eqref{eq:tree-target} at the corresponding generation gives
\eqref{eq:fusion-conclusion} with $b_0=T_0$.
\end{proof}

\begin{proof}[Proof of \Cref{thm:escape}]
Choose a normalized coefficient $\alpha\notin\Lfield$ of $R$
and an order $\nu$ as in \Cref{lem:obstruction}. Let $\tau_R$
be the exponent in \eqref{eq:tau-choice}. These choices depend
only on $R$.

Let $V\subseteq D_R$ be a nonempty open interval. Since $R$ is
nonconstant, $R'$ is a nonzero rational function and has only
finitely many zeros in $D_R$. We may therefore choose a
nondegenerate compact interval $I_*\subset V$ on which
\eqref{eq:derivative-bounds} holds. Applying \Cref{prop:counting}
and then \Cref{prop:fusion} with $f=R$ and $\tau=\tau_R$ yields
a Cantor set $C\subseteq V\cap\Lio$ and a common $b_0$ satisfying
\eqref{eq:escape-intro}.

Fix $\xi\in C$. If $R(\xi)=u/v$ were rational, taking
$(a,b)=(mu,mv)$ with $m$ sufficiently large would contradict
\eqref{eq:escape-intro}. Thus $R(\xi)$ is irrational. For every
$t>\tau_R$, the same bound shows that
$|R(\xi)-a/b|<b^{-t}$ forces $b<b_0$, leaving only finitely
many possible pairs $(a,b)$. Hence $\mu(R(\xi))\leq\tau_R$.
\end{proof}

\section{Classification and consequences}
\label{sec:consequences}

The inclusion $R(\Lio\cap V)\subseteq\Eset$ on a single nonempty
open interval $V\subseteq D_R$ characterizes the coefficient
condition. Excluding rational values then gives the criterion
for Maillet's property on $V$.

\begin{proof}[Proof of \Cref{thm:classification}]
The implication $\textup{(i)}\Rightarrow\textup{(ii)}$ follows
from \Cref{prop:evaluation}. Since $D_R$ contains a nonempty
open interval, $\textup{(ii)}$ implies $\textup{(iii)}$.
Finally, if $R\notin\Lfield(X)$, then \Cref{thm:escape} gives,
in every nonempty open interval $V\subseteq D_R$, a Liouville
number $\xi$ for which $R(\xi)$ is irrational and has finite
irrationality exponent. Thus $R(\xi)\notin\Eset$, proving
$\textup{(iii)}\Rightarrow\textup{(i)}$.
\end{proof}

\begin{proof}[Proof of \Cref{cor:exact}]
Let $U$ be either a nonempty open interval in $D_R$ or $D_R$
itself. Since $\Eset=\Q\cup\Lio$, exact preservation on $U$
is equivalent to the two conditions
\[
 R(\Lio\cap U)\subseteq\Eset,
 \qquad
 R^{-1}(\Q)\cap U\cap\Lio=\varnothing.
\]
By \Cref{thm:classification}, the first is equivalent to
$P,Q\in\Lfield[X]$ in the normalization
\eqref{eq:normalized-pair}. For each $t\in\Q$, coprimality
of $P$ and $Q$ ensures that no zero of $P-tQ$ is a pole of $R$.
Hence the real fiber over $t$ is exactly $Z_{\R}(P-tQ)$, giving
the second condition in \eqref{eq:exact-intro}.
\end{proof}

\subsection{Translations, dilations, and coefficient fields}

\begin{corollary}[Intrinsic descriptions of $\Lfield$]
\label{cor:intrinsic}
For $\alpha\in\R$, the following conditions are equivalent.
\begin{enumerate}[label=\textup{(\roman*)},leftmargin=2.7em]
\item $\alpha\in\Lfield$.
\item $\alpha+\Lio\subseteq\Eset$.
\item $\alpha\Lio\subseteq\Eset$.
\item $\Eset+\alpha=\Eset$.
\end{enumerate}
If $\alpha\ne0$, these conditions are also equivalent to
$\alpha\Eset=\Eset$. Moreover, $\Lfield$ is the largest subfield
$K\subseteq\R$ such that
\begin{equation}\label{eq:universal-extension}
 K(\xi)\subseteq\Eset\qquad\text{for every }\xi\in\Lio.
\end{equation}
\end{corollary}
\begin{proof}
Applying \Cref{thm:classification} to $X+\alpha$ gives
$\textup{(i)}\Leftrightarrow\textup{(ii)}$. Applying it to
$\alpha X$ gives $\textup{(i)}\Leftrightarrow\textup{(iii)}$
when $\alpha\ne0$. For $\alpha=0$, both conditions hold.

Suppose that $\alpha\in\Lfield$. Then
$\alpha+\Q\subseteq\Lfield\subseteq\Eset$, which together with
\textup{(ii)} gives $\alpha+\Eset\subseteq\Eset$. Applying the
same inclusion to $-\alpha$ gives equality. Conversely,
\textup{(iv)} implies \textup{(ii)} since $\Lio\subseteq\Eset$.
For $\alpha\ne0$, the multiplicative assertion follows in the
same way. Indeed, \textup{(i)} and \textup{(iii)} give
$\alpha\Eset\subseteq\Eset$, and applying this inclusion to
$\alpha^{-1}$ gives equality. The reverse implication follows
from \textup{(iii)}.

For $K=\Lfield$, \eqref{eq:universal-extension} follows from
\Cref{prop:evaluation}, since every element of $\Lfield(\xi)$
has the form $R(\xi)$ with $R\in\Lfield(X)$ and $\xi\in D_R$.
Conversely, if $K$ satisfies \eqref{eq:universal-extension},
then every $\alpha\in K$ satisfies
$\alpha\xi\in K(\xi)\subseteq\Eset$ for every $\xi\in\Lio$.
By \textup{(iii)}, $\alpha\in\Lfield$. Thus $K\subseteq\Lfield$.
\end{proof}

\begin{corollary}[The field generated by the coefficients]
\label{cor:coefficient-field}
Let $R\in\R(X)$ be nonconstant, and let $K_R$ be the field
generated over $\Q$ by the coefficients of its normalized pair.
If $R(\Lio\cap V)\subseteq\Eset$ for some nonempty open interval
$V\subseteq D_R$, then $K_R\subseteq\Lfield$. In particular,
every irrational element of $K_R$ is Liouville. For a polynomial
$R(X)=a_0+\cdots+a_mX^m$ satisfying the same hypothesis, this gives
\[
 \Q(a_0,\ldots,a_m)\subseteq\Lfield\subseteq\Eset.
\]
\end{corollary}
\begin{proof}
By \Cref{thm:classification}, every coefficient of the normalized
pair belongs to $\Lfield$. Since $\Lfield$ is a field, it contains
$K_R$. The remaining assertions follow from
\Cref{lem:elementary-field}.
\end{proof}

To see why the distinction between $\Lfield$ and $\Eset$ matters,
choose $\lambda\in\Lio\setminus\Lfield$ by \Cref{prop:size}.
The coefficient field of $X+\lambda$ satisfies
$\Q(\lambda)\subseteq\Lfield(\lambda)\subseteq\Eset$ by
\eqref{eq:universal-extension}. Yet $X+\lambda\notin\Lfield(X)$,
so every nonempty open interval contains a Liouville number $\xi$
for which $\lambda+\xi$ is irrational and has finite irrationality
exponent. The same $\xi$ shows why the proper extension
$K=\Lfield(\lambda)$ fails \eqref{eq:universal-extension}:
$\lambda+\xi\in K(\xi)\setminus\Eset$.

\begin{corollary}[Algebraic coefficients]\label{cor:algebraic}
Let $R\in(\overline\Q\cap\R)(X)$ be nonconstant. The following
conditions are equivalent.
\begin{enumerate}[label=\textup{(\roman*)},leftmargin=2.7em]
\item $R\in\Q(X)$.
\item $R$ has Maillet's property.
\item $R(\Lio\cap V)\subseteq\Eset$ for some nonempty open
      interval $V\subseteq D_R$.
\end{enumerate}
\end{corollary}
\begin{proof}
Maillet's theorem gives $\textup{(i)}\Rightarrow\textup{(ii)}$,
and $\textup{(ii)}\Rightarrow\textup{(iii)}$ is immediate.
Under \textup{(iii)}, the normalized coefficients belong to
$\Lfield$ by \Cref{thm:classification}. They are also algebraic.
Since every irrational element of $\Lfield$ is Liouville and
hence transcendental, all these coefficients are rational.
\end{proof}

\begin{remark}[A projective formulation]
Put $\Eset_\infty=\Eset\cup\{\infty\}\subseteq\mathbb P^1(\R)$.
For a nonconstant real rational map acting on $\mathbb P^1(\R)$,
\[
 R(\Eset_\infty)\subseteq\Eset_\infty
 \quad\Longleftrightarrow\quad R\in\Lfield(X).
\]
Indeed, a map defined over $\Lfield$ sends $\mathbb P^1(\Q)$
into $\mathbb P^1(\Lfield)\subseteq\Eset_\infty$. At a Liouville
number its value is either $\infty$ at a pole or an element of
$\Eset$ by \Cref{prop:evaluation}. Conversely, projective
preservation gives $R(\Lio\cap D_R)\subseteq\Eset$, so
\Cref{thm:classification} applies. For the induced map $R:\mathbb P^1(\R)\to\mathbb P^1(\R)$, we have
\[
 R(\Lio)\subseteq\Lio
 \quad\Longleftrightarrow\quad
 R\in\Lfield(X)
 \ \text{and}\
 R^{-1}\bigl(\mathbb P^1(\Q)\bigr)\cap\Lio=\varnothing.
\]
The condition at $\infty$ excludes poles in $\Lio$.
\end{remark}

\subsection{Degree one and rational fibers}

\begin{corollary}[M\"obius rigidity]\label{cor:mobius}
A real M\"obius transformation $M$ has Maillet's property on
some nonempty open interval in $D_M$ if and only if $M\in\Q(X)$.
If $M\in\Lfield(X)\setminus\Q(X)$, then
\[
 \{\xi\in\Lio\cap D_M:M(\xi)\in\Q\}
\]
is dense in $D_M$.
\end{corollary}
\begin{proof}
Suppose first that $M\in\Lfield(X)\setminus\Q(X)$, and fix a
nonempty open interval $V\subseteq D_M$. The restriction of $M$
to $V$ is a homeomorphism onto the open interval $M(V)$. For
every $t\in\Q\cap M(V)$, its unique preimage belongs to
$\Lfield$, since $M^{-1}$ is also defined over $\Lfield$.

At most two such preimages can be rational. Indeed, a M\"obius
transformation taking three distinct rational points to rational
points belongs to $\Q(X)$. To see this, send each ordered triple
to $(0,1,\infty)$ by a M\"obius transformation over $\Q$ and use
uniqueness of the map between the two triples. Thus, with at most
two exceptions, the preimages of $\Q\cap M(V)$ belong to
$\Lfield\setminus\Q\subseteq\Lio$. Since $\Q\cap M(V)$ remains
dense after the removal of two points, these Liouville preimages
are dense in $V$.

Now suppose that $M$ has Maillet's property on a nonempty open
interval. By \Cref{thm:classification}, $M\in\Lfield(X)$.
The density assertion rules out $M\notin\Q(X)$, and the converse
follows from Maillet's theorem.
\end{proof}

In particular, if $M\notin\Q(X)$, every nonempty open interval
$V\subseteq D_M$ contains a Liouville number $\xi$ for which
$M(\xi)\notin\Lio$: this follows from \Cref{thm:escape}
when $M\notin\Lfield(X)$ and from the density of rational
fibers above otherwise. In higher degree, a point in a rational
fiber need only be algebraic over the coefficient field.
Even the fiber of $X^2$ over $2$ contains points outside $\Lfield$.

\subsection{Real-analytic and entire functions}

The local rigidity theorem of Part~I reduces the real-analytic
case to the rational classification.

\begin{corollary}[Real-analytic functions]\label{cor:analytic}
Let $U\subseteq\R$ be a nonempty open interval, and let
$f:U\to\R$ be nonconstant and real-analytic. Then
$f(\Lio\cap U)\subseteq\Eset$ if and only if $f=R|_U$ for some
$R\in\Lfield(X)$ with $U\subseteq D_R$.

Moreover, $f(\Lio\cap U)\subseteq\Lio$ if and only if such a
representation exists and its normalized pair $R=P/Q$ satisfies
\[
 Z_{\R}(P-tQ)\cap U\cap\Lio=\varnothing
 \qquad(t\in\Q).
\]
For $f\equiv c$, the two inclusions are equivalent to
$c\in\Eset$ and $c\in\Lio$, respectively.
\end{corollary}
\begin{proof}
Suppose that $f$ has no rational representation on $U$.
By \cite[Theorem~1.1]{MarquesI}, there exist $\xi\in U\cap\Lio$,
a finite $\tau>2$, and an integer $b_0$ such that
\[
 \abs{f(\xi)-\frac ab}>b^{-\tau}
 \qquad(a\in\Z,\ b\in\Z,\ b\geq b_0).
\]
This excludes a rational value of $f(\xi)$, since any rational
number admits representations with arbitrarily large denominators.
It also gives $\mu(f(\xi))\leq\tau$. Thus $f(\xi)\notin\Eset$.

Consequently, $f(\Lio\cap U)\subseteq\Eset$ forces $f=R|_U$
for a real rational function $R$ without poles on $U$.
For such an $R$, \Cref{thm:classification} gives the first
equivalence and \Cref{cor:exact} gives the second. The constant
cases follow since $U\cap\Lio$ is nonempty.
\end{proof}

\begin{corollary}[Entire functions]\label{cor:entire}
Let $F:\C\to\C$ be nonconstant and entire. Then
\begin{equation}\label{eq:entire-weak}
 F(\Lio)\subseteq\Eset
 \quad\Longleftrightarrow\quad F\in\Lfield[X].
\end{equation}
Moreover,
\[
 F(\Lio)\subseteq\Lio
 \quad\Longleftrightarrow\quad
 F\in\Lfield[X]\ \text{and}\ F^{-1}(\Q)\cap\Lio=\varnothing.
\]
\end{corollary}
\begin{proof}
Suppose that $F(\Lio)\subseteq\Eset$. Since $\Eset\subseteq\R$
and $\Lio$ is dense in $\R$, continuity gives $F(\R)\subseteq\R$.
The restriction $F|_{\R}$ is nonconstant by the identity theorem.
Applying \Cref{cor:analytic} on $\R$, we obtain $F=P/Q$ there,
with $P,Q\in\Lfield[X]$ coprime and $Q$ monic.

The entire function $QF-P$ vanishes on $\R$, so it vanishes
on $\C$. Any complex zero of $Q$ would therefore be a zero of
$P$, contrary to coprimality by B\'ezout's identity.
Hence $Q=1$ and $F\in\Lfield[X]$. Conversely, every polynomial
over $\Lfield$ maps $\Lio$ into $\Eset$ by \Cref{prop:evaluation}.
Exact preservation is then equivalent to
$F^{-1}(\Q)\cap\Lio=\varnothing$.
\end{proof}

\begin{corollary}[Holomorphic functions on connected domains]
\label{cor:holomorphic}
Let $\Omega\subseteq\C$ be a connected open set with
$\Omega\cap\R\ne\varnothing$, and let
$F:\Omega\to\C$ be nonconstant and holomorphic. Then
$F(\Lio\cap\Omega)\subseteq\Eset$ if and only if there is a
rational function $R\in\Lfield(X)$ with no pole in $\Omega$
such that $F=R|_\Omega$.
Moreover, $F(\Lio\cap\Omega)\subseteq\Lio$ if and only if
such an $R$ exists and
\[
 R^{-1}(\Q)\cap\Omega\cap\Lio=\varnothing.
\]
For constant functions $F\equiv c$, the two conditions are
respectively $c\in\Eset$ and $c\in\Lio$.
\end{corollary}
\begin{proof}
Suppose first that $F(\Lio\cap\Omega)\subseteq\Eset$.
Choose a nonempty real open interval $I\subseteq\Omega\cap\R$.
Because $\Lio\cap I$ is dense in $I$, continuity and
$\Eset\subseteq\R$ give $F(I)\subseteq\R$.
The restriction $F|_I$ is real-analytic and nonconstant,
by the identity theorem. By \Cref{cor:analytic},
$F=R$ on $I$, where $R=P/Q\in\Lfield(X)$ is reduced
and $Q$ is monic. The holomorphic function $QF-P$
vanishes on $I$ and hence identically on $\Omega$.
If $Q$ vanished at any $z\in\Omega$, the identity would
give $P(z)=0$, contradicting the coprimality of $P,Q$.
Thus $R$ has no poles in $\Omega$ and $F=R|_\Omega$.
Conversely, if $R\in\Lfield(X)$ has no pole in $\Omega$,
\Cref{prop:evaluation} shows that
$R(\Lio\cap\Omega)\subseteq\Eset$.
The exact-preservation criterion is obtained by excluding
rational values. The constant case is immediate.
\end{proof}

\section{Strong Liouville coefficients and rational $U_m$-fibers}
\label{sec:strong}

We now prove \Cref{thm:strong-nonlinear}. The construction uses
continued fractions whose denominators have prescribed prime
valuations. These valuations force irreducibility in every degree,
while the growth of consecutive denominators yields strong
Liouville approximation.

Recall that a real irrational number $\sigma$ is strong Liouville
if the denominators $q_n$ of its continued-fraction convergents
satisfy $\log q_{n+1}/\log q_n\to\infty$; see
\eqref{eq:strong}. For an algebraic number $\beta\in\C$, let
$\mathrm h(\beta)$ denote its absolute logarithmic Weil height.
For $m\geq1$, a complex transcendental number $\kappa$ belongs
to $U_m$ if $m$ is the least positive integer for which there
exist algebraic numbers $\beta_n$ of degree exactly $m$ and
positive real numbers $w_n$ such that
\[
 \mathrm h(\beta_n)\longrightarrow\infty,\qquad
 w_n\longrightarrow\infty,\qquad
 |\kappa-\beta_n|\leq e^{-w_n\mathrm h(\beta_n)}.
\]
This is LeVeque's absolute-logarithmic-height formulation;
see \cite[\S~1]{BiluMarques2026} and
\cite[Remark~2.6]{BiluMarquesMoreira2026}.
It agrees with the usual Mahler--Koksma classification
\cite[Theorem~3.6]{Bugeaud2004}. In particular,
$U_1=\Lio$ and $U_m\cap U_1=\varnothing$ for $m\geq2$.

We use the following consequence of the gap principle in
\cite[Theorem~3.1 and Corollary~3.3]{BiluMarques2026}.

\begin{lemma}[$U_m$ gap criterion]\label{lem:Um-gap}
Let $m\geq1$, let $\kappa\in\C$ be transcendental, and let
$(\beta_n)$ be algebraic numbers of degree at most $m$.
Set $H_n=\mathrm h(\beta_n)$, and suppose that there are
positive real numbers $w_n$ and constants $B,C>0$ such that
\[
 H_n\longrightarrow\infty,\qquad w_n\longrightarrow\infty,
 \qquad |\kappa-\beta_n|\leq Ce^{-w_nH_n},
 \qquad H_{n+1}\leq Bw_nH_n.
\]
If $[\Q(\beta_n):\Q]=m$ for all sufficiently large $n$,
then $\kappa\in U_m$.
\end{lemma}

This is precisely the combination of
\cite[Theorem~3.1 and Corollary~3.3]{BiluMarques2026};
no monotonicity or pairwise distinctness assumption on the
approximants is needed in addition to the stated hypotheses.

\begin{proof}[Proof of \Cref{thm:strong-nonlinear}]
Enumerate the nonzero rationals as
\[
\begin{gathered}
 \Q^\times=\{r_1,r_2,\ldots\},\qquad
 r_j=\frac{A_j}{B_j},\\
 A_j\in\Z\setminus\{0\},\qquad B_j\in\Z_{>0},\qquad
 \gcd(A_j,B_j)=1.
\end{gathered}
\]
We construct a continued fraction
\[
 \sigma=[1;a_1,a_2,\ldots],\qquad
 \frac{p_n}{q_n}=[1;a_1,\ldots,a_n],
\]
with $q_{-1}=0$ and $q_0=1$. Suppose that
$a_1,\ldots,a_{n-1}$ have been chosen. Select a prime
$\ell_n$ such that
\begin{equation}\label{eq:strong-prime-choice}
 \ell_n\nmid q_{n-1}\prod_{j=1}^n A_jB_j.
\end{equation}
Since $q_{n-1}$ is invertible modulo $\ell_n^2$, the
recurrence $q_n=a_nq_{n-1}+q_{n-2}$ allows us to choose
$a_n$ arbitrarily large subject to
\begin{equation}\label{eq:strong-prime-valuation}
 q_n\equiv\ell_n\pmod{\ell_n^2},\qquad
 q_n>q_{n-1}^{n+1}.
\end{equation}
In particular, $v_{\ell_n}(q_n)=1$, while
$\ell_n\nmid p_n$ because $\gcd(p_n,q_n)=1$.
The growth condition in \eqref{eq:strong-prime-valuation}
shows that $\log q_{n+1}/\log q_n\to\infty$.
Thus $\sigma$ is strong Liouville and $1<\sigma<2$.

Fix $j\geq1$ and $m\geq2$. For every $n\geq j$,
\begin{equation}\label{eq:strong-Eisenstein}
 G_{j,m,n}(X):=B_jp_nX^m-A_jq_n
\end{equation}
is Eisenstein at $\ell_n$: by
\eqref{eq:strong-prime-choice} and \eqref{eq:strong-prime-valuation},
the prime $\ell_n$ divides the constant coefficient
exactly once and does not divide the leading coefficient.
Hence every root of
\[
 X^m-\frac{r_j}{p_n/q_n}
\]
has degree exactly $m$ over $\Q$.

Let $\beta\in\C$ be any root of $\beta^m=r_j/\sigma$.
The number $\beta$ is transcendental since $\sigma$ is.
The implicit-function theorem, applied to $zw^m=r_j$
near $(\sigma,\beta)$, supplies a holomorphic branch
$u(z)$ with $u(\sigma)=\beta$ and $u(z)^m=r_j/z$.
For all sufficiently large $n$, put
$\beta_n=u(p_n/q_n)$. By the preceding Eisenstein
argument, $[\Q(\beta_n):\Q]=m$, and the standard
continued-fraction estimate gives
\begin{equation}\label{eq:strong-root-approximation}
 |\beta-\beta_n|
 \ll_{j,m,\beta}\left|\sigma-\frac{p_n}{q_n}\right|
 <\frac{C_{j,m,\beta}}{q_nq_{n+1}}.
\end{equation}

Set $H_n=\mathrm h(\beta_n)$. The identity
$\beta_n^m=A_jq_n/(B_jp_n)$ implies
\[
 mH_n=\mathrm h\left(\frac{A_jq_n}{B_jp_n}\right).
\]
Since $1<p_n/q_n<2$ for all sufficiently large $n$ and
the greatest common divisor
of $A_jq_n$ and $B_jp_n$ divides $|A_jB_j|$, we have
\begin{equation}\label{eq:strong-root-height}
 H_n=\frac1m\log q_n+O_{j,m}(1).
\end{equation}
In particular, $H_n\to\infty$. Define
\[
 w_n:=\frac{\log q_{n+1}}{2H_n}.
\]
The strong Liouville condition and
\eqref{eq:strong-root-height} imply $w_n\to\infty$.
Equation \eqref{eq:strong-root-approximation} yields
\[
 |\beta-\beta_n|\leq C_{j,m,\beta}e^{-w_nH_n}.
\]
Furthermore, since $m\geq2$,
\eqref{eq:strong-root-height} gives, for all sufficiently
large $n$,
\[
 H_{n+1}\leq\log q_{n+1}=2w_nH_n.
\]
The gap criterion in \Cref{lem:Um-gap} now shows that
$\beta\in U_m$. As $j,m$, and the chosen root $\beta$
were arbitrary, every nonzero rational fiber of
$\sigma X^m$ consists of $U_m$-numbers.

By \Cref{prop:strong}, $\sigma\in\Lfield$.
For $m\geq2$ the polynomial $\sigma X^m$ therefore
belongs to $\Lfield[X]$, and its rational fiber over zero
is $\{0\}$. The other rational fibers are disjoint from
$\Lio=U_1$ by what we have just proved. Applying
\Cref{cor:exact} gives
\[
 \sigma\xi^m\in\Lio\qquad
 (\xi\in\Lio,\ m\geq2).
\]

Finally, at each step of the continued-fraction construction,
the required congruence class modulo $\ell_n^2$ contains
infinitely many admissible positive integers $a_n$. Choosing
between two at every step gives continuum many distinct
continued fractions, all of which have the stated properties.
\end{proof}

This also provides a coefficient-field interpretation of the
strong Liouville multiplication phenomenon in
\cite{Petruska1992}. By \Cref{prop:fixed-gap-inclusion}, all the
coefficients constructed above belong to $\Ggap$. The next
section shows that preservation in higher degree also occurs
with coefficients outside the entire field $\Q(\Ggap)$.

\section{An explicit coefficient beyond fixed-gap approximation}
\label{sec:examples}

We construct a positive number $\alpha$ outside the field generated
by the fixed-gap class defined below such that $\alpha X^m$ has
Maillet's property for every $m\geq2$. The growth of the
denominators places $\alpha$ in $\Lfield$ and separates it from
that field, while their prime factors allow us to control the
rational fibers of each monomial.

\subsection{The fixed-gap class}

For a reduced fraction $a/b$ with $b>0$, write
$\height(a/b)=\log\max\{|a|,b\}$. Let $\Ggap$ consist of the
real numbers $\lambda$ for which there exist rationals $u_j$,
positive reals $w_j$, and a constant $A>0$ satisfying
\begin{equation}\label{eq:fixed-gap}
 \begin{gathered}
 \height(u_j)\longrightarrow\infty,\qquad
 w_j\longrightarrow\infty,\qquad
 |\lambda-u_j|\leq e^{-w_j\height(u_j)},\\
 \height(u_{j+1})\leq A w_j\height(u_j).
 \end{gathered}
\end{equation}
This is the class denoted by $\mathcal L$ in
\cite[Section~1.1]{BiluMarques2026}, following
\cite{ChavesMarquesTrojovsky2021}. A single constant $A$ controls
the growth of the heights along the entire sequence.

\begin{proposition}\label{prop:fixed-gap-inclusion}
The class $\Ggap$ contains every strong Liouville number and
satisfies $\Ggap\subseteq\Lfield\setminus\Q$.
\end{proposition}
\begin{proof}
Let $\sigma$ be strong Liouville, with convergents $u_j=p_j/q_j$.
Write $h_j=\height(u_j)=\log q_j+O(1)$. After discarding finitely
many terms, set $w_j=(\log q_{j+1})/h_j$. By \eqref{eq:strong},
$w_j\to\infty$. The convergent estimate gives
$|\sigma-u_j|<q_{j+1}^{-1}=e^{-w_jh_j}$, and
$h_{j+1}\leq2\log q_{j+1}=2w_jh_j$ for all sufficiently large
$j$. Thus $\sigma\in\Ggap$.

Now let $\lambda\in\Ggap$, write $h_j=\height(u_j)$ for a
sequence satisfying \eqref{eq:fixed-gap}, and fix an integer
$N\geq1$. Discard an initial segment so that $w_j\geq2N$.
For every sufficiently large $H$, let $j$ be the first remaining
index with $h_j>\log H$. Its predecessor $i=j-1$ then satisfies
$h_i\leq\log H$.

If $w_i h_i\geq2N\log H$, then $\den(u_i)\leq H$ and
\[
 |\lambda-u_i|\leq H^{-2N}\leq(H\den(u_i))^{-N}.
\]
Otherwise, $h_j\leq Aw_i h_i<2AN\log H$, so
$\den(u_j)\leq H^{2AN}$. Since $e^{h_j}>H$ and
$\den(u_j)\leq e^{h_j}$, we have
\[
 |\lambda-u_j|\leq e^{-2Nh_j}\leq(H\den(u_j))^{-N}.
\]
Thus $D=\lceil\max\{1,2AN\}\rceil$ is admissible in
\Cref{def:field}.

Finally, if $\lambda=a/b$ were rational, then $u_j\ne\lambda$
for all sufficiently large $j$, since $h_j\to\infty$. Rational
separation would give
$|\lambda-u_j|\geq(b\den(u_j))^{-1}\geq b^{-1}e^{-h_j}$,
contrary to \eqref{eq:fixed-gap} and $w_j\to\infty$.
\end{proof}

\subsection{Construction of the coefficient}

Let $\ell_n$ be the $n$th prime and let $v_2(n)$ denote the
exponent of $2$ in $n$. Define
\begin{equation}\label{eq:example-recurrence}
 e_n:=3+v_2(n),\qquad Q_1:=2,\qquad
 Q_{n+1}:=\ell_{n+1}Q_n^{e_n},
\end{equation}
and put
\begin{equation}\label{eq:example-alpha}
 \alpha:=1+\sum_{n\geq1}\frac1{Q_n},\qquad
 u_n:=1+\sum_{j=1}^n\frac1{Q_j}=\frac{P_n}{Q_n}.
\end{equation}

\begin{theorem}\label{thm:example}
The number $\alpha$ in \eqref{eq:example-alpha} belongs to
$\Lfield\setminus\Q(\Ggap)$ and satisfies
\[
 \alpha\xi^m\in\Lio
 \qquad(\xi\in\Lio,\ m\in\Z_{\geq2}).
\]
Moreover, for each $m\geq2$, the irrationality exponents of the
real roots of $\alpha X^m=t$, as $t$ ranges over $\Q^\times$,
are bounded by a constant depending only on $m$.
\end{theorem}

\begin{lemma}\label{lem:example-basic}
The fractions $u_n=P_n/Q_n$ are reduced and satisfy
\begin{equation}\label{eq:prime-valuation}
 v_{\ell_n}(Q_n)=1,\qquad \ell_n\nmid P_n.
\end{equation}
Moreover, $1<u_n<\alpha<2$ and
\begin{equation}\label{eq:example-tail}
 \frac1{Q_{n+1}}<\alpha-u_n<\frac2{Q_{n+1}}.
\end{equation}
Writing $t_n=\log Q_n$, we have
\begin{equation}\label{eq:log-growth}
 e_nt_n\leq t_{n+1}\leq(e_n+2)t_n.
\end{equation}
In particular, $\alpha$ is Liouville.
\end{lemma}
\begin{proof}
The numerators satisfy $P_1=3$ and
\[
 P_{n+1}=P_n\ell_{n+1}Q_n^{e_n-1}+1.
\]
Thus $P_{n+1}$ is congruent to $1$ modulo every prime divisor
of $Q_{n+1}$. The prime divisors of $Q_n$ are precisely
$\ell_1,\ldots,\ell_n$, with $\ell_n$ occurring to the first
power. This proves reducedness and \eqref{eq:prime-valuation}.

Since $Q_{n+1}\geq Q_n^3\geq4Q_n$, comparison with a geometric
series gives $\alpha<2$ and \eqref{eq:example-tail}. Every prime
divisor of $Q_n+1$ is at least $\ell_{n+1}$, so
$\ell_{n+1}\leq Q_n+1\leq Q_n^2$. Taking logarithms in
\eqref{eq:example-recurrence} gives \eqref{eq:log-growth}.

Finally, $e_{2^j}=j+3$. Along these indices,
\eqref{eq:example-tail} gives
$0<\alpha-u_n<2Q_n^{-e_n}$, and hence rational approximations
of arbitrarily high order. Thus $\alpha\in\Lio$.
\end{proof}

\begin{lemma}[Approximation at a prescribed scale]
\label{lem:example-scale}
For each integer $N\geq2$, there is an integer $D_N\geq1$ such
that every sufficiently large real $H$ admits an index $n$ with
\begin{equation}\label{eq:example-scale}
 Q_n\leq H^{D_N},\qquad Q_{n+1}\geq(HQ_n)^N.
\end{equation}
Consequently, $\alpha\in\Lfield$.
\end{lemma}
\begin{proof}
Set $s_N=2N$ and $g_N=2^{s_N-3}$. The indices $i$
for which $e_i\geq s_N$ are exactly the multiples of $g_N$.
If $i$ and $i'=i+g_N$ are consecutive such indices, then
$e_j\leq s_N-1$ for $i<j<i'$. Iterating
\eqref{eq:log-growth} over these $g_N-1$ indices gives
\begin{equation}\label{eq:between-good-indices}
 t_{i'}\leq C_Nt_{i+1},\qquad C_N:=(s_N+1)^{g_N-1}.
\end{equation}

For sufficiently large $H$, choose consecutive multiples $i,i'$
of $g_N$ with $Q_i\leq H<Q_{i'}$.
If $Q_{i+1}\geq H^{2N}$, take $n=i$. Then $Q_n\leq H$ and
$Q_{n+1}\geq H^{2N}\geq(HQ_n)^N$.

If $Q_{i+1}<H^{2N}$, take $n=i'$. By
\eqref{eq:between-good-indices},
\[
 Q_n\leq Q_{i+1}^{C_N}<H^{2NC_N}.
\]
Since $Q_n>H$ and $e_n\geq2N$, we also have
$Q_{n+1}\geq Q_n^{2N}>(HQ_n)^N$. Thus $D_N=2NC_N$
works in both cases.

For a prescribed order $M\geq1$ in \Cref{def:field}, apply
\eqref{eq:example-scale} with $N=M+1$. By \eqref{eq:example-tail},
the resulting truncation satisfies
\[
 Q_n\leq H^{D_{M+1}},\qquad
 |\alpha-u_n|<2(HQ_n)^{-(M+1)}\leq(HQ_n)^{-M}
\]
for all sufficiently large $H$. Hence $\alpha\in\Lfield$.
\end{proof}

\begin{lemma}\label{lem:not-fixed-gap}
The number $\alpha$ lies outside the field $\Q(\Ggap)$.
\end{lemma}
\begin{proof}
We first show that every reduced fraction $a/b$ with $b\geq2$,
other than the truncations $u_n$, satisfies
\begin{equation}\label{eq:nontruncation-lower}
 \left|\alpha-\frac ab\right|\geq\frac1{8b^3}.
\end{equation}
Choose $n$ so that $Q_n\leq b<Q_{n+1}$. If
$b\leq Q_{n+1}/(4Q_n)$, separation from $u_n$ and
\eqref{eq:example-tail} give
\[
 \left|\alpha-\frac ab\right|
 >\frac1{bQ_n}-\frac2{Q_{n+1}}
 \geq\frac1{2bQ_n}\geq\frac1{2b^2}.
\]
Otherwise, compare $a/b$ with $u_{n+1}$. Since $b<Q_{n+1}$ and
$Q_{n+2}\geq Q_{n+1}^3\geq4bQ_{n+1}$, we obtain
\[
 \left|\alpha-\frac ab\right|
 >\frac1{bQ_{n+1}}-\frac2{Q_{n+2}}
 \geq\frac1{2bQ_{n+1}}>\frac1{8b^3},
\]
where the last inequality follows from
$Q_{n+1}<4bQ_n\leq4b^2$. This proves
\eqref{eq:nontruncation-lower}.

For $\gamma\in\Lfield$, let $d_\gamma(N)$ be the least admissible
denominator exponent in \eqref{eq:scale-property} at order $N$.
Let $\mathcal K_{\mathrm{pol}}$ be the set of those $\gamma$
for which $d_\gamma(N)$ is bounded by a polynomial in $N$.
The proof of \Cref{prop:fixed-gap-inclusion} gives
$d_\gamma(N)=O_\gamma(N)$ for $\gamma\in\Ggap$, so
$\Ggap\subseteq\mathcal K_{\mathrm{pol}}$.

We claim that $\mathcal K_{\mathrm{pol}}$ is a field. It contains $\Q$, since exact
approximations are allowed. For $\beta,\gamma\in K$, the proof
of \Cref{lem:synchronization} gives a common denominator exponent
$D_1+D_2+D_1D_2$ at order $M$, where
$D_1=d_\beta(M)$ and $D_2=d_\gamma(M(D_1+1))$.
This exponent is polynomially bounded in $M$. The proof of
\Cref{prop:field} uses $M=2N+1$ for subtraction and multiplication,
at most doubling the common exponent. For inversion it uses
$d_\beta(N+1)+1$ when $\beta\ne0$. Thus $\mathcal K_{\mathrm{pol}}$ is a field and
$\Q(\Ggap)\subseteq\mathcal K_{\mathrm{pol}}$.

It remains to show that $\alpha\notin\mathcal K_{\mathrm{pol}}$.
We prove that
\begin{equation}\label{eq:example-required-exponent}
 d_\alpha(N)\geq3^{\,2^{N-4}-1}\qquad(N\geq5).
\end{equation}
Fix such an $N$, put $g_N=2^{N-4}$, and take $H=Q_{i+1}$
with $i$ an arbitrarily large multiple of $g_N$. For all sufficiently
large such $i$, choose a fraction satisfying
\eqref{eq:scale-property} with $D=d_\alpha(N)$, and reduce it
to $a/b$. Reduction preserves both bounds. These fractions
converge to the irrational number $\alpha$, so $b\to\infty$.
Since $(Hb)^{-N}<1/(8b^3)$ for all sufficiently large $i$,
\eqref{eq:nontruncation-lower} forces $a/b=u_n$ for some $n$.
By \eqref{eq:example-tail},
\[
 Q_{n+1}>(HQ_n)^N.
\]
In particular, $n>i$. By \eqref{eq:log-growth}, $e_n+2>N$.
As $e_n=3+v_2(n)$, this gives $g_N\mid n$, hence $n\geq i+g_N$.
Since $Q_{k+1}\geq Q_k^3$ and $u_n$ is reduced,
\[
 H^{d_\alpha(N)}\geq b=Q_n\geq Q_{i+g_N}
 \geq Q_{i+1}^{\,3^{g_N-1}}=H^{\,3^{g_N-1}}.
\]
This proves \eqref{eq:example-required-exponent}. Its lower bound
grows faster than every polynomial in $N$, so
$\alpha\notin\mathcal K_{\mathrm{pol}}$ and therefore
$\alpha\notin\Q(\Ggap)$.
\end{proof}

\subsection{The rational fibers of $\alpha X^m$}

\begin{lemma}\label{lem:example-fibers}
Let $m\geq2$ be an integer and let $t\in\Q^\times$. Every real
root $\beta$ of $\alpha X^m=t$ satisfies
\begin{equation}\label{eq:fiber-exponent}
 \mu(\beta)\leq m+D_m,
\end{equation}
where $D_m$ is the exponent in \Cref{lem:example-scale}.
\end{lemma}
\begin{proof}
Write $t=A/B$ with $A\in\Z\setminus\{0\}$ and $B\in\Z_{>0}$,
and consider
\[
 G_n(X):=BP_nX^m-AQ_n.
\]
For all sufficiently large $n$, the prime $\ell_n$ divides
neither $A$ nor $B$. By \eqref{eq:prime-valuation},
\[
 v_{\ell_n}\left(\frac{AQ_n}{BP_n}\right)=1.
\]
Thus $AQ_n/(BP_n)$ cannot be an $m$th power in $\Q$, and $G_n$
has no rational root. In particular, $b^mG_n(a/b)$ is a nonzero
integer for every $a\in\Z$ and $b\in\Z_{>0}$, so
\begin{equation}\label{eq:integer-polynomial-lower}
 |G_n(a/b)|\geq b^{-m}.
\end{equation}

For each sufficiently large $b$, choose $n=n(b)$ by
\Cref{lem:example-scale} with $H=b$ and $N=m$. Then
\[
 Q_n\leq b^{D_m},\qquad Q_{n+1}\geq(bQ_n)^m.
\]
The second inequality forces $n(b)\to\infty$, so
\eqref{eq:integer-polynomial-lower} applies for all sufficiently
large $b$. Since $\beta^m=A/(B\alpha)$, we have
\[
 G_n(\beta)=\frac{AQ_n}{\alpha}(u_n-\alpha).
\]
Using \eqref{eq:example-tail} and the choice of $n$, we obtain
\[
 |G_n(\beta)|
 <\frac{2|A|}{\alpha}\,b^{-m}Q_n^{1-m}
 =o(b^{-m}).
\]
Here $m\geq2$ and $Q_n\to\infty$.

On $[\beta-1,\beta+1]$, the bound $P_n<2Q_n$ gives
$|G_n'|\ll_{m,t,\beta}Q_n$. Thus, for every sufficiently large
$b$ and every $a$ with $|\beta-a/b|\leq1$, the mean value
theorem yields
\[
 Q_n|\beta-a/b|
 \gg_{m,t,\beta}|G_n(a/b)-G_n(\beta)|
 \geq\tfrac12 b^{-m}.
\]
Since $Q_n\leq b^{D_m}$, it follows that
$|\beta-a/b|\gg_{m,t,\beta} b^{-(m+D_m)}$. The same bound holds
when $|\beta-a/b|>1$. This proves \eqref{eq:fiber-exponent}.
\end{proof}

\begin{proof}[Proof of \Cref{thm:example}]
By \Cref{lem:example-basic,lem:example-scale,lem:not-fixed-gap},
we have $\alpha\in\Lfield\setminus\Q(\Ggap)$. Fix $m\geq2$.
The fiber of $\alpha X^m$ over zero is $\{0\}$, and
\Cref{lem:example-fibers} excludes Liouville numbers from every
nonzero rational fiber. Since $\alpha X^m\in\Lfield[X]$,
\Cref{cor:exact} gives $\alpha\xi^m\in\Lio$ for every
$\xi\in\Lio$. The bound $m+D_m$ in \eqref{eq:fiber-exponent}
is independent of $t$.
\end{proof}

Together with \Cref{prop:field,prop:fixed-gap-inclusion,prop:size},
this gives
\[
 \Q(\Ggap)\subsetneq\Lfield\subsetneq\Eset.
\]

\section{Scope of the characterization}
\label{sec:scope}

By \cite[Corollary~1.2]{MarquesI}, every nonconstant real-analytic
function with Maillet's property on an open interval $U$ is
the restriction of a real rational function with no pole on $U$.
For a nonconstant $R\in\R(X)$, \Cref{thm:classification} shows
that the inclusion $R(\Lio\cap U)\subseteq\Eset$ on a single
nonempty open interval $U\subseteq D_R$ forces $R\in\Lfield(X)$.
The same inclusion then holds throughout $D_R$. Exact
preservation on $U$ is equivalent to the additional condition
$R^{-1}(\Q)\cap U\cap\Lio=\varnothing$.

The exponent in \Cref{thm:escape} depends on $R$ through the
bidegree and the obstruction order in \eqref{eq:tau-choice}.
Part~I gives an absolute exponent for nonrational real-analytic
functions \cite[Theorem~1.1]{MarquesI}. For rational functions,
dependence on the coefficients is necessary even in degree one.

\begin{proposition}\label{prop:no-uniform-escape}
For every integer $N\geq2$ there exists
$\beta\in\R\setminus\Lfield$ such that $\xi+\beta$ is irrational
and $\mu(\xi+\beta)\geq N$ for every $\xi\in\Lio$.
\end{proposition}
\begin{proof}
Fix $N\geq2$. Define $\beta=[0;2,a_2,\ldots]$ recursively,
with initial convergent denominators $q_0=1$, $q_1=2$,
by setting $a_{j+1}=q_j^{2N-1}$ for $j\geq1$.
The denominator recurrence gives
\[
 q_j^{2N}\leq q_{j+1}=q_j^{2N}+q_{j-1}\leq2q_j^{2N}.
\]
The continued-fraction formula
$\mu(\beta)=1+\limsup_{j\to\infty}\log q_{j+1}/\log q_j$
therefore gives $\mu(\beta)=2N+1$. Hence $\beta\notin\Lfield$
by \Cref{lem:elementary-field}.

Fix $\xi\in\Lio$ and put $z=\xi+\beta$. If $z$ were rational,
$\beta=z-\xi$ would be Liouville. Let $W=N(2N+1)+1$, and choose
reduced fractions $p/b$ with $b\to\infty$ and
$|\xi-p/b|<b^{-W}$. Let $a/h$ be the first convergent of $\beta$
with $h\geq b$. Then $b\leq h\leq2b^{2N}$ for all sufficiently
large $b$, and $|\beta-a/h|<h^{-(2N+1)}$.
For $s_b=p/b+a/h$, we obtain
\[
 0<|z-s_b|\leq b^{-W}+h^{-(2N+1)}
          \leq2(bh)^{-N}
          \leq2\den(s_b)^{-N}
\]
for all sufficiently large $b$. Since $s_b\to z$ and $z$ is
irrational, the reduced denominators tend to infinity.
For every $t<N$, the inequality
$0<|z-s_b|<\den(s_b)^{-t}$ therefore holds for all sufficiently
large $b$. Thus $\mu(z)\geq N$.
\end{proof}

For $R=P/Q\in\Lfield(X)$, the criterion in \Cref{cor:exact}
requires the exclusion of Liouville roots of $P-tQ$ for every
$t\in\Q$. Applying this structural criterion to a particular
map requires arithmetic information about its rational fibers.

In higher degree, \Cref{thm:strong-nonlinear} supplies
continuum many positive strong Liouville coefficients $\sigma$
for which $\sigma X^m$ has Maillet's property for every
$m\geq2$, with each nonzero rational fiber consisting of
$U_m$-numbers. These coefficients belong to $\Ggap$.
The explicit coefficient in \Cref{thm:example} lies instead in
$\Lfield\setminus\Q(\Ggap)$, yet every monomial
$\alpha X^m$, $m\geq2$, still has Maillet's property.
Its nonzero rational fibers satisfy the stronger condition of
having uniformly bounded irrationality exponents for each $m$.
The linear map $\alpha X$ takes the rational value $1$ at
$\alpha^{-1}\in\Lfield\setminus\Q\subseteq\Lio$.

\subsection*{Further questions}
\noindent
The criterion in \Cref{cor:exact} is structural rather than
algorithmic. Even for the family $R(X)=\lambda X^m$, it is natural
to ask for intrinsic conditions on $\lambda\in\Lfield^\times$
which ensure that no nonzero rational fiber contains a Liouville
number. Equivalently, one seeks to characterize the coefficients
$\lambda$ for which
\[
 \lambda\xi^m\in\Lio
 \qquad(\xi\in\Lio).
\]
The examples of \Cref{thm:strong-nonlinear,thm:example}
show that such coefficients exist in very different approximation
classes, while the choice $\lambda=\mu^{-m}$,
$\mu\in\Lfield\setminus\Q$, shows that not every
nonrational coefficient in $\Lfield$ works.
A related problem is to identify natural probability measures on
structured subclasses of $\Lfield$ under which the prevalence
of the fiber condition can be studied. Ordinary Lebesgue
almost-everywhere statements on $\Lfield$ would be vacuous,
since $\dim_{\mathrm H}\Lfield=0$.

\section*{Use of artificial intelligence tools}

OpenAI's ChatGPT was used for English-language editing and organization
of the exposition. The author is responsible for
the correctness of all statements and proofs.

\end{document}